\documentclass[preprint]{elsarticle}
\makeatletter
\def\ps@pprintTitle{%
 \let\@oddhead\@empty
 \let\@evenhead\@empty
 \def\@oddfoot{}%
 \let\@evenfoot\@oddfoot}
\makeatother

\usepackage{amsmath}
\usepackage{amssymb}
\usepackage{amsthm}
\usepackage{graphicx}
\usepackage{epic,eepic,epsfig}
\usepackage{color}
\usepackage{placeins}	
\usepackage{multirow}
\usepackage{epstopdf}
\usepackage{varwidth}
\usepackage[table]{xcolor}
\usepackage{hyperref}
\usepackage{subcaption}

\usepackage{tikz}
\usetikzlibrary{decorations.pathreplacing}

\newcommand{\rr}[1]{\textcolor{black} {#1}}

\definecolor{tabclr}{cmyk}{0,0,1,0}

\begin{document}

\title{Learning adaptive coarse spaces of BDDC algorithms
for stochastic elliptic problems with oscillatory and high contrast coefficients}

\author[CUHK]{Eric Chung\corref{cor}\fnref{fn1}}
\ead{tschung@math.cuhk.edu.hk}

\author[KHU]{Hyea Hyun Kim\fnref{fn2}}
\ead{hhkim@khu.ac.kr}


\author[CUHK]{Ming Fai Lam}
\ead{mflam@math.cuhk.edu.hk}

\author[CUHK]{Lina Zhao}
\ead{lzhao@math.cuhk.edu.hk}

\cortext[cor]{Corresponding author}
\fntext[fn1]{The research of Eric Chung is supported by Hong Kong RGC General Research Fund (Project number 14317516)
and CUHK Faculty of Science Research Incentive Fund 2015-16.}
\fntext[fn2]{The research of Hyea Hyun Kim is
supported by the National Research Foundation of Korea(NRF) grants
funded by NRF-\rr{2015R1A5A1009350} and by \rr{NRF-2019R1A2C1010090}.}

\address[CUHK]{Department of Mathematics, The Chinese University of Hong Kong, Hong Kong SAR}
\address[KHU]{Department of Applied Mathematics and Institute of Natural Sciences, Kyung Hee University,
Korea}

\begin{abstract}
In this paper, we consider the balancing domain decomposition by constraints (BDDC) algorithm with adaptive coarse spaces for
a class of stochastic elliptic problems. The key ingredient in the construction of the coarse space is the solutions of
local spectral problems, which depend on the coefficient of the PDE. This poses a significant challenge for stochastic coefficients
as it is computationally expensive to solve the local spectral problems for every realisation of the coefficient.
To tackle this computational burden, we propose a machine learning approach. Our method is based on the use of a deep neural network (DNN)
to approximate the relation between the stochastic coefficients and the coarse spaces.
For the input of the DNN, we apply the Karhunen-Lo\`eve expansion and use the first few dominant terms in the expansion. The output
of the DNN is the resulting coarse space, which is then applied with the standard adaptive BDDC algorithm.
We will present some numerical results with oscillatory and high contrast coefficients to show the efficiency and robustness of the proposed scheme.
\end{abstract}

\begin{keyword}
BDDC, stochastic partial differential equation, artificial neural network, coarse space, high contrast
\end{keyword}

\maketitle

\section{Introduction}

Coefficient functions are often one of the main difficulties in modeling a real life problem. Due to the nature of the phenomena, these coefficients can be oscillatory and with high contrast features, all such conditions hinder the development of a fast and accurate numerical solution. Under some conditions, numerical methods like discontinuous Galerkin methods could give robust scheme as in \cite{chung2018adaptive,chung2014adaptive}. However, more importantly, coefficients are not constant \rr{and they usually involve} uncertainties and randomness. Therefore, a stochastic partial differential equation (SPDE) is considered for the modeling of these stochastic coefficient functions. There are several popular approaches for solving SPDE numerically, for example, Monte Carlo simulation, stochastic Galerkin method, \rr{and} stochastic collocation method in \cite{babuvska2007stochastic,babuska2004galerkin,ghanem2003stochastic}. All these methods require a high computational power in realistic simulations and cannot be generalized to similar situations. In the past decades, machine learning technique has been applied on many aspects such as image classification, fluid dynamics and solving differential equations in \cite{brunton2020machine,kutz2017deep,lu2021deepxde,zhao2016spectral,heinlein2021combining,vasilyeva2020learning,wang2020deep,chung2021multi,yeung2020deep}. In this paper, instead of directly \rr{solving} the SPDE by machine learning, we choose a moderate step that \rr{involves} the advantage of machine learning and also the accuracy from traditional numerical solver.

We first apply the Karhunen-Lo\`eve expansion on the stochastic coefficient as in \cite{dostert2006coarse,wheeler2011multiscale,zhang2004efficient} to capture the characteristic with finite terms. After this expansion, the SPDE is reduced to a deterministic differential equation and can be solved with a balancing domain decomposition by constraints (BDDC) algorithm, which can  prevent the computational problem caused by the refinement of spatial mesh resolution and give an accurate solution rather than just \rr{using} machine learning technique without suitably designed conditions. After obtaining sufficient number of random solutions from \rr{the} BDDC method, we can derive the statistics of the original SPDE solution. However, the cost for forming and solving generalized eigenvalue problems in the BDDC algorithm with good accuracy is considerable especially for high dimensional problems. We note that such eigenvalue problems are considered to enrich the coarse
space adaptivity in the BDDC algorithm.
Therefore, we introduce neural network with \rr{a} user-defined number \rr{of layers and} neurons to compensate the increase in computational cost with unexpected times of BDDC solutions.

The adaptive BDDC algorithm with enriched primal unknowns is considered in this paper as its ability on oscillatory and  high contrast coefficients has been shown in \cite{kim2017bddc,kim2015bddc}. Among different domain decomposition methods, the considered adaptive BDDC does not require a strong assumption on coefficient functions and the subdomain partitions to achieve a good performance like the standard BDDC algorithm \cite{BDDC-D-03,BDDC-Mandel-Dohrmann-Tezaur}. It is because the additional coarse space basis functions computed by the dominant eigenfunctions are related to the ill-conditioning in the standard BDDC algorithm with highly varying coefficients. The introduction of these dominant eigenfunctions \rr{thus} could greatly improve the robustness of numerical scheme for problems with rough and high contrast coefficients. Moreover, an estimate of condition number of this adaptive BDDC algorithm can be shown to be only controlled by a given tolerance without any extra assumptions on the coefficients or subdomain partitions. For simplicity, in the remaining of \rr{the} paper, we will call the resulting new numerical scheme combining the machine learning technique and adaptive BDDC algorithm as learning adaptive BDDC algorithm.

Finally, to end this section, we state the main idea of our proposed scheme and the model problem. In this paper, we integrate an artificial neural network with learning abilities into a BDDC algorithm with adaptively enriched coarse spaces for efficient solutions of the neural network approximation of stochastic elliptic problems. Let $\mathcal{D}$ be a bounded domain in $\mathbb{R}^d$, $d=2$, and $\Omega$ be a stochastic probability space:
\begin{equation}\label{model:pb}
\begin{split}
-\nabla \cdot ( \rho(\boldsymbol{x},\omega) \nabla u) &= f, \quad \text{ in } \mathcal{D}, \\
u &= 0, \quad \text{ on } \partial \mathcal{D},
\end{split}
\end{equation}
where $\boldsymbol{x}\in \mathcal{D},$ $\omega\in\Omega$ and $\rho(\boldsymbol{x},\omega)$ is uniformly positive and is highly heterogeneous with very high contrast. Here, we only assume the two-dimensional spatial domain, however, the three-dimensional case is also applicable and is shown to be robust for deterministic elliptic problems in \cite{kim2017bddc}. Moreover, as the model problem \eqref{model:pb} can be extended to fluid flow problems, we seldom call the coefficient function  $\rho(\boldsymbol{x},\omega)$ as permeability function.

The rest of the paper is organized as follows. In Section~\ref{sec:BDDC}, a brief formulation of \rr{the} BDDC algorithm with adaptively enriched coarse problems is presented for two-dimensional elliptic problems. Then, the Karhunen-Lo\`eve expansion is introduced and the details of the artificial neural network used \rr{are} clarified in Section~\ref{sec:ML}. Two concrete examples with an explicit analytical expression of the Karhunen-Lo\`eve expansion, some network training and testing parameters are presented in the first half of Section~\ref{sec:num} and results of various numerical experiments are
presented subsequently. Finally, a concluding remark is given in Section~\ref{sec:conclusion}.

\section{\label{sec:BDDC}Adaptive BDDC algorithm}

The main feature of \rr{the} adaptive BDDC is the local generalized eigenvalue problem defined on every subdomain interface, which introduces the adaptively enriched primal unknowns. In this section, we repeat the formulation presented by the first and second authors in \cite{kim2017bddc} and give a brief overview of the adaptive BDDC algorithm. We refer to \cite{kim2017bddc,kim2015bddc} for further details, where the analysis of \rr{the} condition number estimation and the robustness of numerical scheme with oscillatory and high contrast coefficients can be found.

\subsection{Local linear system}

We first introduce a discrete form of the model problem~\eqref{model:pb} in a deterministic
fashion. Let $V_h$ be the space of conforming linear finite element functions
with respect to a given mesh on $\mathcal{D}$ with \rr{the} mesh size $h>0$ and with the zero value on $\partial \mathcal{D}$.
We will then find the approximate solution $u\in V_h$ such that
\begin{equation}
a(u,v) = (f,v), \quad\quad \forall v\in V_h,
\label{eq:cg}
\end{equation}
where
\begin{equation}
a(u,v) = \int_{\mathcal{D}} \rho(\boldsymbol{x}) \nabla u \cdot \nabla v \, dx, \quad
(f,v) = \int_{\mathcal{D}} f \, v \, dx.
\end{equation}
We assume that the spatial domain $\mathcal{D}$ is partitioned into a set of $N$ non-overlapping subdomains $\{ \mathcal{D}_i \}$, $i=1,2,\cdots, N$,
so that $\mathcal{D} = \cup_{i=1}^N \mathcal{D}_i$.
We note that the subdomain boundaries do not cut triangles equipped for $V_h$.
We allow the coefficient $\rho(\boldsymbol{x})$ to have high contrast jumps and oscillations across subdomains
and on subdomain interfaces.
Let $a_i(u,v)$ be the bilinear form of the model elliptic problem (\ref{eq:cg}) restricted to each subdomain $\mathcal{D}_i$ defined as
$$a_i(u,v)=\int_{\mathcal{D}_i} \rho(\boldsymbol{x}) \nabla u \cdot \nabla v \, dx, \quad \forall u, v \in X_i,$$
where $X_i$ is the restriction of $V_h$ to $\mathcal{D}_i$.

In the BDDC algorithm, the original problem~\eqref{eq:cg} is reduced to a subdomain interface problem and solved by an iterative method combined with a preconditioner. The interface problem can be obtained by solving a Dirichlet problem in each subdomain. After choosing dual and primal unknowns on the subdomain
interface unknowns, the interface problem is then solved by utilizing local problems and one global
coarse problem corresponding to the chosen sets of dual and primal unknowns, respectively.
At each iteration, the residuals are multiplied by certain scaling factors to balance
the errors across the subdomain interface regarding to the energy of each subdomain problem.
The coarse problem aims to correct the global part of the error in each iteration and thus the choice of primal unknowns is important in obtaining a good performance as the number of subdomains increases.
The basis for primal unknowns is obtained by the minimum energy extension for a given constraint
at the location of primal unknowns and such a basis provides a robust coarse problem with a good
energy estimate.
We refer to \cite{BDDC-D-03,LW-FETIDP-BDDC,BDDC-Mandel-Dohrmann-Tezaur,TW-Book} for general
introductions to BDDC algorithm.

\subsection{\label{BDDC:FETI:theory}Notation and preliminary results}
To facilitate our discussion, we first introduce some notation.
Let $S^{(i)}$ be the Schur complement matrix obtained from the local stiffness matrix $A^{(i)}$
after eliminating unknowns interior to $\mathcal{D}_i$, where $A^{(i)}$ is defined by
$a_i(u,v)=u^T A^{(i)} v$, for all $u,v\in X_i$.
In the following we will use the same symbol to represent a finite element function and its corresponding coefficient vector
in order to simplify the notations.

Recall that $X_i$ is the restriction of the finite element space $V_h$ to each subdomain $\mathcal{D}_i$.
Let $W_i$ be the restriction of $X_i$ to $\partial \mathcal{D}_i$.
We then introduce the product spaces
$$X=\prod_{i=1}^N X_i,\quad W=\prod_{i=1}^N W_i,$$
where we remark that the functions in $X$ and $W$ are totally decoupled across the subdomain interfaces.
In addition, we introduce partially coupled subspaces $\widetilde{X}$, $\widetilde{W}$,
and fully coupled subspaces $\widehat{X}$, $\widehat{W}$,
where some primal unknowns are strongly coupled for functions in $\widetilde{X}$ or $\widetilde{W}$,
while the functions in $\widehat{X}$, $\widehat{W}$ are fully coupled across the subdomain interfaces.

Next, we present basic description of the BDDC algorithm; see~\cite{BDDC-D-03,LW-FETIDP-BDDC,BDDC-Mandel-Dohrmann-Tezaur,TW-Book}.
For simplicity, the two-dimensional case will be considered.
After eliminating unknowns interior to each subdomain, the Schur complement matrices $S^{(i)}$
are obtained from $A^{(i)}$ and they form the algebraic problem considered in the BDDC algorithm,
which is to find $\widehat{w} \in \widehat{W}$ such that
\begin{equation}\label{primal-pb}\sum_{i=1}^N R_i^T S^{(i)} R_i \widehat{w}= \sum_{i=1}^N R_i^T g_i,
\end{equation}
where $R_i: \widehat{W} \rightarrow W_i$ is the restriction operator into $\partial \mathcal{D}_i$,
and $g_i \in W_i$ depends on the source term $f$.

The BDDC preconditioner is built based on the partially coupled space $\widetilde{W}$.
Let $\widetilde{R}_i: \widetilde{W} \rightarrow W_i$ be the restriction into $\partial \mathcal{D}_i$
and let $\widetilde{S}$ be the partially coupled matrix defined by
$$\widetilde{S}=\sum_{i=1}^N \widetilde{R}_i^T S^{(i)} \widetilde{R}_i.$$
For the space $\widetilde{W}$, we can express it as the product of the two spaces
$$\widetilde{W}=W_{\Delta} \times \widehat{W}_{\Pi},$$
where $\widehat{W}_{\Pi}$ consists of vectors of the primal unknowns and
$W_{\Delta}$ consists of vectors of dual unknowns, which are strongly coupled
at the primal unknowns and decoupled at the remaining interface unknowns, respectively.
We define $\widetilde{R}: \widehat{W} \rightarrow \widetilde{W}$ such that
$$\widetilde{R}=\begin{pmatrix}
R_{\Delta} \\
R_{\Pi}
\end{pmatrix},$$
where $R_{\Delta}$ is the mapping from $\widehat{W}$ to $W_{\Delta}$ and $R_{\Pi}$ is the restriction
from $\widehat{W}$ to $\widehat{W}_{\Pi}$.
We note that $R_{\Delta}$ is obtained as
$$R_{\Delta}=\begin{pmatrix} R_{\Delta}^{(1)} \\ R_{\Delta}^{(2)}\\ \vdots \\ R_{\Delta}^{(N)}
\end{pmatrix},$$
where $R_{\Delta}^{(i)}$ is the restriction from $\widehat{W}$ to $W_{\Delta}^{(i)}$
and $W_{\Delta}^{(i)}$ is the space of dual unknowns of $\mathcal{D}_i$.

The BDDC preconditioner is then given by
\begin{equation}\label{precond:BDDC}
M^{-1}_{BDDC}=\widetilde{R}^T \widetilde{D} \widetilde{S}^{-1} \widetilde{D}^T \widetilde{R},
\end{equation}
where $\widetilde{D}$ is a scaling matrix of the form
$$\widetilde{D}=\sum_{i=1}^N \widetilde{R}_i^T D_i \widetilde{R}_i.$$
Here the matrices  $D_i$ are defined for unknowns in $W_i$ and they
are introduced to make the preconditioner robust to the heterogeneity in $\rho(\boldsymbol{x})$ across the subdomain interface.
In more detail, $D_i$ consists of blocks $D_F^{(i)}$ and $D_V^{(i)}$, where $F$ denotes an equivalence class
shared by two subdomains, i.e., $\mathcal{D}_i$ and its neighboring subdomain $\mathcal{D}_j$, and
$V$ denotes the end points of $F$, respectively.
We call such equivalence classes $F$ and $V$ as edge and vertex in two dimensions, respectively. 
In our BDDC algorithm, unknowns at subdomain vertices are included to the set of primal unknowns
and adaptively selected primal constraints are later included to the set after a change of basis formulation.
For a given edge $F$ in two dimensions, the matrices $D_F^{(l)}$ and $D_V^{(l)}$
satisfy a partition of unity property, i.e., $D_F^{(i)}+D_F^{(j)}=I$ and $\sum_{l \in n(V) } D_V^{(l)}=1$,
where $n(V)$ denotes the set of subdomain indices sharing the vertex $V$. The matrices $D_F^{(l)}$ and $D_V^{(l)}$
are called scaling matrices.
As mentioned earlier, the scaling matrices help to balance the residual error
at each iteration with respect to the energy of subdomain problems sharing the interface.
For the case when $\rho(\boldsymbol{x})$ is identical across the interface $F$, $D_F^{(i)}$ and $D_F^{(j)}$ are chosen
simply as multiplicity scalings, i.e., $1/2$, but for a general case when $\rho(\boldsymbol{x})$ has discontinuities,
different choice of scalings, such as $\rho$-scalings or deluxe scalings, can be more effective.
The scaling matrices $D_V^{(l)}$ can be chosen using similar ideas.
We refer to \cite{KRR-2015} and references therein for scaling matrices.

We note that by using the definitions of $\widetilde{R}$ and $\widetilde{S}$, the matrix in the left hand side of \eqref{primal-pb}
can be written as
$$\sum_{i=1}^N R_i^T S^{(i)} R_i = \widetilde{R}^T \widetilde{S} \widetilde{R}.$$
In the BDDC algorithm, the system in \eqref{primal-pb} is solved by an iterative method
with the preconditioner~\eqref{precond:BDDC}. Thus its performance is analyzed by estimating
the condition number of
\begin{equation}\label{cond:BDDC}M^{-1}_{BDDC} \widetilde{R}^T \widetilde{S} \widetilde{R}=\widetilde{R}^T \widetilde{D} \widetilde{S}^{-1} \widetilde{D}^T \widetilde{R}\widetilde{R}^T \widetilde{S} \widetilde{R}.
\end{equation}

\subsection{Generalized eigenvalue problems}
After the preliminaries, we are ready to state the generalized eigenvalue problems which introduce the adaptive enrichment of coarse components.
For  an equivalence class $F$ shared by two subdomains
$\mathcal{D}_i$ and $\mathcal{D}_j$, the following generalized eigenvalue problem is proposed in \cite{Klawonn:PAMM:2014}:
\begin{equation}\label{Klawonn:EIG}
\Big((D_F^{(j)})^TS_F^{(i)}D_F^{(j)} + (D_F^{(i)})^T S_F^{(j)} D_F^{(i)} \Big) v
=\lambda \, \Big( \widetilde{S}_F^{(i)}:\widetilde{S}_F^{(j)} \Big) v,
\end{equation}
where $S_F^{(i)}$ and $D_F^{(i)}$ are the block matrices of $S^{(i)}$ and $D_i$ corresponding to unknowns
interior to $F$, respectively.
The matrix $\widetilde{S}_F^{(i)}$ are the Schur complement
of $S^{(i)}$ obtained after eliminating unknowns except those interior to $F$. In addition, for symmetric and semi-positive definite matrices $A$ and $B$,
their parallel sum $A:B$ is defined by, see \cite{Parallel:Sum},
\begin{equation}
A:B=B(A+B)^{+}A, \label{eq:psum}
\end{equation}
where $(A+B)^{+}$ is a pseudo inverse of $A+B$.
We note that the problem in \eqref{Klawonn:EIG} is identical
to that considered in \cite{Pechstein:slide:2013} when $D_F^{(i)}$ are chosen as
the deluxe scalings, i.e.,
$$D_F^{(i)}=(S_F^{(i)}+S_F^{(j)})^{-1} S_F^{(i)}.$$
\rr{We note that a similar generalized eigenvalue problem is considered and extended to three-dimensional problems in the first and second authors' work \cite{kim2017bddc}.}

\section{\label{sec:ML}Learning adaptive BDDC algorithm}

In addition to the  BDDC algorithm, we  perform the Karhunen-Lo\`eve (KL) expansion to decompose the stochastic permeability $\rho(\boldsymbol{x},\omega)$ for the preparation of the learning adaptive BDDC algorithm. As shown in \cite{zhang2004efficient}, only a few terms in the KL expansion are already enough to approximate the stochastic permeability with a reasonable accuracy, which can be efficiently reduce computational cost. Once we obtain the KL series expansion, we can easily produce a certain amount of sample data as a training set for the artificial neural network. After the training process, the neural network captures the main characteristics of the dominant eigenvectors from the adaptive BDDC algorithm, and neural network approximation can be obtained. In our proposed algorithm, we still need some \rr{resulting} samples from the deterministic BDDC algorithm to train the neural network, which is a process of a relatively high computational cost. However, after this setup, we can apply this network to obtain \rr{an} approximate solution directly without the use of other stochastic sampling methods such as Monte Carlo simulation, which involves even  higher computational cost for the generation of large number of random samples or realizations.

In this section, we discuss the proposed algorithm in details. Moreover, our resulting neural network can  be \rr{also} applied on problems with similar stochastic properties. It \rr{thus} saves the cost of training again and we can readily use the trained neural network for prediction. Some testing samples will be generated to test the network performance \rr{and} to show its generalization capacity in the later section.

\subsection{Karhunen-Lo\`eve expansion}

To start our learning adaptive algorithm, we first apply the KL expansion on the stochastic field $\rho(\boldsymbol{x},\omega)$. To ensure the positivity of stochastic permeability field, a logarithmic transformation is considered as $K(\boldsymbol{x},\omega) = \log(\rho(\boldsymbol{x},\omega))$. Let $C_K(\boldsymbol{x},\boldsymbol{\hat{x}})$ be the covariance function of $K(\boldsymbol{x},\omega)$ at two locations $\boldsymbol{x}$ and $\boldsymbol{\hat{x}}$. Since the covariance function $C_K(\boldsymbol{x},\boldsymbol{\hat{x}})$ is symmetric and positive definite, it can be decomposed
into:
\begin{align}
C_K(\boldsymbol{x},\boldsymbol{\hat{x}}) = \sum_{i=1}^\infty \lambda_i f_i(\boldsymbol{x}) f_i(\boldsymbol{\hat{x}}),
\end{align}
where $\lambda_i$ and $f_i$ are eigenvalues and eigenfunctions computed from the following Fredholm integral equation:
\begin{align}
  \int_\mathcal{D} C_K(\boldsymbol{x},\boldsymbol{\hat{x}})f(\boldsymbol{x})\;d\boldsymbol{x} = \lambda f(\boldsymbol{\hat{x}}),
\end{align}
where $\{f_i(\boldsymbol{x})\}$ are orthogonal and deterministic functions. It is noted that for some specially chosen covariance functions, we can find the $(\lambda,f)$ eigenpair analytically. We will list some examples in the \rr{section of} numerical results. After solving the Fredholm integral equation, we are ready to express the KL expansion of the log permeability:
\begin{align}\label{KLexpre}
K(\boldsymbol{x},\omega) = E[K](\boldsymbol{x})+\sum_{i=1}^\infty\sqrt{\lambda_i}\xi_i f_i(\boldsymbol{x}).
\end{align}
Here $E[K](\boldsymbol{x})$ is the expected value of $K(\boldsymbol{x},\omega)$ and ${\xi_i}$ are identically independent distributed Gaussian random variables with mean 0 and variance 1. As mentioned above, a few terms in KL expansion can give a reasonably accurate approximation, which can be explained now as we can just pick the dominant eigenfunctions in \eqref{KLexpre} and the other eigenvalues will decay rapidly.

Therefore, KL expansion is an efficient method to represent the stochastic coefficient and flavour the subsequent process. Although analytical  eigenfunctions and eigenvalues cannot always be found, there are various numerical methods to compute the KL coefficients in \cite{schwab2006karhunen,wang2008karhunen}. After the generation of some sets of $\{\xi_i\}$, which are the only stochastic variables in the permeability after KL expansion, as the input of neural network, we plug back KL expanded permeability coefficients into the adaptive BDDC algorithm, and the resulting dominant eigenfunctions in the coarse spaces will be the target output of the our designated neural network.

\subsection{Neural network}

\begin{figure}[h]
	\centering
	\begin{tikzpicture}[shorten >=1pt]
		\tikzstyle{unit}=[draw,shape=circle,minimum size=1.15cm]
		\tikzstyle{hidden}=[draw,shape=circle,minimum size=1.15cm]

		\node[unit](x0) at (0,3.5){$x_1$};
		\node[unit](x1) at (0,2){$x_2$};
		\node at (0,1){\vdots};
		\node[unit](xd) at (0,0){$x_{R}$};

		\node[hidden](h10) at (3,4){$y_1^{(1)}$};
		\node[hidden](h11) at (3,2.5){$y_2^{(1)}$};
		\node at (3,1){\vdots};
		\node[hidden](h1m) at (3,-0.5){$y_{n^{(1)}}^{(1)}$};

		\node[hidden](h22) at (5,-0.5){$y_{n^{(2)}}^{(2)}$};
		\node at (5,1){\vdots};
		\node[hidden](h21) at (5,2.5){$y_2^{(2)}$};
		\node[hidden](h20) at (5,4){$y_1^{(2)}$};
		
		\node(d3) at (6,-0.5){$\ldots$};
		\node(d2) at (6,2.5){$\ldots$};
		\node(d1) at (6,4){$\ldots$};

		\node[hidden](hL12) at (7,-0.5){\small{$y_{n^{(L-1)}}^{(L-1)}$}};
		\node at (7,1){\vdots};
		\node[hidden](hL11) at (7,2.5){\small{$y_2^{(L-1)}$}};
		\node[hidden](hL10) at (7,4){\small{$y_1^{(L-1)}$}};
		
		\node[hidden](hL0) at (9,4){$y_1^{(L)}$};
		\node[hidden](hL1) at (9,2.5){$y_2^{(L)}$};
		\node at (9,1){\vdots};
		\node[hidden](hLm) at (9,-0.5){$y_{n^{(L)}}^{(L)}$};

		\node[unit](y1) at (12,3.5){$y_1^{(L+1)}$};
		\node[unit](y2) at (12,2){$y_2^{(L+1)}$};
		\node at (12,1){\vdots};	
		\node[unit](yc) at (12,0){$y_O^{(L+1)}$};

		\draw[->] (x0) -- (h11);
		\draw[->] (x0) -- (h1m);
        \draw[->] (x0) -- (h10);

		\draw[->] (x1) -- (h11);
		\draw[->] (x1) -- (h1m);
        \draw[->] (x1) -- (h10);

		\draw[->] (xd) -- (h11);
		\draw[->] (xd) -- (h1m);
        \draw[->] (xd) -- (h10);

		\draw[->] (hL0) -- (y1);
		\draw[->] (hL0) -- (yc);
		\draw[->] (hL0) -- (y2);

		\draw[->] (hL1) -- (y1);
		\draw[->] (hL1) -- (yc);
		\draw[->] (hL1) -- (y2);

		\draw[->] (hLm) -- (y1);
		\draw[->] (hLm) -- (y2);
		\draw[->] (hLm) -- (yc);

		\draw[->] (h10) -- (h21);
		\draw[->] (h10) -- (h22);
		\draw[->] (h10) -- (h20);

		
		\draw[->] (h11) -- (h21);
		\draw[->] (h11) -- (h22);
		\draw[->] (h11) -- (h20);
		
		\draw[->] (h1m) -- (h21);
		\draw[->] (h1m) -- (h22);
		\draw[->] (h1m) -- (h20);

		\draw[->] (hL10) -- (hL0);
		\draw[->] (hL11) -- (hL0);
		\draw[->] (hL12) -- (hL0);
		
		\draw[->] (hL10) -- (hL1);
		\draw[->] (hL11) -- (hL1);
		\draw[->] (hL12) -- (hL1);
		
		\draw[->] (hL10) -- (hLm);
		\draw[->] (hL11) -- (hLm);
		\draw[->] (hL12) -- (hLm);
		
		\draw [decorate,decoration={brace,amplitude=10pt},xshift=-4pt,yshift=0pt] (-0.5,4.1) -- (0.75,4.1) node [black,midway,yshift=+0.6cm]{input layer};
		\draw [decorate,decoration={brace,amplitude=10pt},xshift=-4pt,yshift=0pt] (2.5,4.5) -- (3.75,4.5) node [black,midway,yshift=+0.6cm]{$1^{\text{st}}$ hidden layer};
		\draw [decorate,decoration={brace,amplitude=10pt,mirror},xshift=-4pt,yshift=0pt] (8.5,-1) -- (9.75,-1) node [black,midway,yshift=-0.6cm]{$L^{\text{th}}$ hidden layer};
        \draw [decorate,decoration={brace,amplitude=10pt,mirror},xshift=-4pt,yshift=0pt] (4.5,-1) -- (5.75,-1) node [black,midway,yshift=-0.6cm]{$2^{\text{nd}}$ hidden layer};
		\draw [decorate,decoration={brace,amplitude=10pt},xshift=-4pt,yshift=0pt] (6.5,4.5) -- (7.75,4.5) node [black,midway,yshift=+0.6cm]{$(L-1)^{\text{th}}$ hidden layer};
		\draw [decorate,decoration={brace,amplitude=10pt},xshift=-4pt,yshift=0pt] (11.5,4.1) -- (12.75,4.1) node [black,midway,yshift=+0.6cm]{output layer};
	\end{tikzpicture}
	\caption[Network graph for a $(L+1)$-layer perceptron.]{An illustration of neural network structure of a $(L+1)$-layer perceptron with $R$ input neurons and $O$ output neurons. The $l^{\text{th}}$ hidden layer contains $n^{(l)}$ hidden neurons.}
	\label{fig:FFNN}
\end{figure}
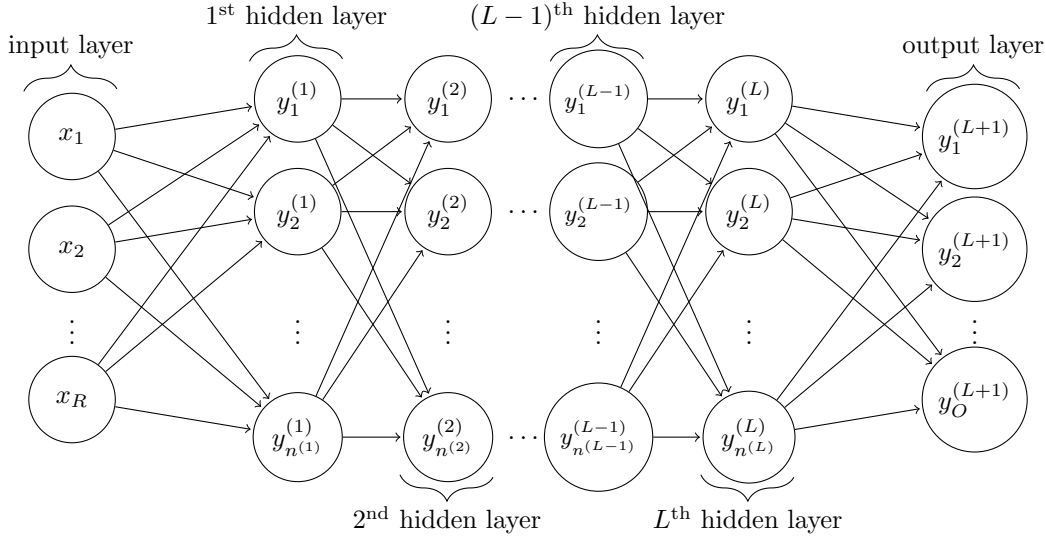

There are many types of neural networks for different usages. Since in our learning adaptive BDDC algorithm the desired neural network is to capture  numeric feature from the data as a supervised learning, \rr{ a fully connected feedforward neural network is chosen.} An illustration of the neural network structure is shown in Figure \ref{fig:FFNN}. In a general feedforward neural network, there are $L$ hidden layers between the input and output layer, and the arrows from layers to layers indicate that all the neurons in the previous layer are used to produce every neuron in the next layer. \rr{Computational cost increases with $L$, but there are no significant improvements in predictions of our tests, $L=1$ is thus chosen throughout numerical experiments.} 

We denote the number of neurons in the input layer to be $R$, which is the number of terms in the truncated KL expansion. The inputs are actually the gaussian random variables $\{\xi_i\}$ in one realization. The output of the network is a column vector \rr{that} consists of all dominant \rr{eigenvectors} in the coarse space. The number of neurons in the output layer is denoted as $O$, which is the sum of length of all resulting eigenvectors obtained from the deterministic BDDC algorithm. This value $O$ depends on many factors, such as the geometry of the spatial domain, the structure of the grid, and parameters of the BDDC algorithm.

In the training, we use \rr{the} scaled conjugate gradient \rr{algorithm in \cite{moller1993scaled}} to do \rr{the} task of minimization. Therefore, unlike the  gradient descent method, we do not need to determine the value of learning rate at each step. However, the initial choice \rr{for the} scaled conjugate gradient \rr{algorithm} may still affect the performance of \rr{training}. Before training, we have to generate training samples and decide some conditions on stopping criteria. When the minimization process is completed for all the training samples, we are said to finish one epoch. There are two stopping criteria that are usually considered. \rr{They} are the minimum \rr{for the} cost function gradient and the maximum number of epochs trained. When either of these conditions are satisfied, the training is stopped and the network will be tested for performance on a testing set.

In both training and testing results,  in order to obtain an error measure less \rr{affected by scales of data set}, a normalized root mean squared error (NRMSE) is used as an error of reference to see how well \rr{are} the network estimation and prediction:
\begin{align*}
&RMSE = \sqrt{\dfrac{\sum_{i=1}^M\frac{1}{O}||y(\boldsymbol{\xi}^{(i)})-f_{NN}(\boldsymbol{\xi}^{(i)})||_2^2}{M}}, \\
&NRMSE = \dfrac{RMSE}{\max_{i=1,2,\dots,M}||y(\boldsymbol{\xi}^{(i)})-f_{NN}(\boldsymbol{\xi}^{(i)})||_\infty},
\end{align*}
where $\boldsymbol{\xi}^{(i)}\sim N(0,1)\in\mathbb{R}^R$ is a column vector of random variables \rr{that} follow the standard normal distribution in the $i$-th realization, $y(\cdot)$ is the dominant eigenvectors obtained from the adaptive BDDC algorithm when the KL expanded stochastic permeability function is used, and $f_{NN}(\cdot)$ is the regression function describing the neural network (NN). Here, $O$ is the dimension of network output, and $M$ is the number of samples. 
For simplicity, the  considered NRMSE can be understood as the root mean squared error divided by the maximum element of the absolute difference of all the eigenvectors.

Overall, the proposed scheme can be concluded as follows:
\begin{itemize}
  \item Step 1: Perform KL expansion on the  logarithmic stochastic permeability function $K(\boldsymbol{x},\omega)$
  \item Step 2: Generate $M$ realizations  of \rr{$\{\boldsymbol{\xi}^{(i)}\}$} and obtain the corresponding  BDDC dominant eigenvectors \rr{$\{y(\boldsymbol{\xi}^{(i)})\}$}, which are the training data for the neural network
  \item Step 3: Define training conditions and train the neural network
  \item Step 4: Examine the network performance whether the NRMSE is satisfied, otherwise, go back step 3 and change training conditions
\end{itemize}

\section{\label{sec:num}Numerical results}

In this section, supporting numerical results are presented to show the performance of our proposed learning adaptive BDDC algorithm. We will consider various choices of permeability coefficients $\rho(\boldsymbol{x},\omega)$ with two major stochastic behaviours. As the concern of our numerical tests are on the learning algorithm, we fix all the parameters that are not related. In all the experiments, $\mathcal{D}\in\mathbb{R}^2$ is chosen to be a unit square spatial domain, and it is partitioned into 16 uniform square subdomains with $H=4^{-1}$ as the coarse grid size. Each subdomain is then further divided into uniform grids with \rr{a} fine grid size $h=32^{-1}$. In the following, we will first describe the considered stochastic coefficients and training conditions.

\subsection{Choices of stochastic coefficients}

From the KL expansion expression \eqref{KLexpre}, we know the eigenvalues and eigenfunctions computed from the Fredholm integral equation \rr{are related} to the covariance function of $K(\boldsymbol{x},\omega)$, which can be treated as the starting point of KL expansion. Throughout the experiments, two specially chosen covariance functions whose eigenvalues and eigenfunctions in the Fredholm integral equation can be computed analytically are considered with different expected values. \rr{For} $\boldsymbol{x}=(x_1,x_2)^T$ \rr{and} $\boldsymbol{\hat{x}}=(\hat{x}_1,\hat{x}_2)^T$, these two covariance functions are
\begin{enumerate}
  \item Brownian sheet covariance function:
  $$
  C_K(\boldsymbol{x},\boldsymbol{\hat{x}}) = \min(x_1,\hat{x}_1)\min(x_2,\hat{x}_2),
  $$
  \rr{where} the corresponding  eigenvalues and eigenfunctions are
  \begin{align*}
    \lambda_k &= \dfrac{16}{((2i-1)^2\pi^2)((2j-1)^2\pi^2)},\\
    f_k(\boldsymbol{x}) &= 2\sin\left((i-\dfrac{1}{2})\pi x_1\right)\sin\left((j-\dfrac{1}{2})\pi x_2\right).
  \end{align*}

  \item Exponential covariance function:
  $$
  C_K(\boldsymbol{x},\boldsymbol{\hat{x}}) = \sigma_K^2\exp\left(-\dfrac{|x_1-\hat{x}_1|}{\eta_1}-\dfrac{|x_2-\hat{x}_2|}{\eta_2}\right),
  $$
  where $\sigma_K^2$ and $\eta_i$ are the variance and  correlation length in the $x_i$ direction of the process, respectively. The corresponding  eigenvalues and eigenfunctions are
  \begin{align*}
    \lambda_k &= \dfrac{4\eta_1\eta_2\sigma_K^2}{(r_{1,i}^2\eta_1^2+1)(r_{2,j}^2\eta_2^2+1)},\\
    f_k(\boldsymbol{x}) &= \dfrac{r_{1,i} \eta_1 \cos(r_{1,i}x_1)+\sin(r_{1,i}x_1)}{\sqrt{(r_{1,i}^2\eta_1^2+1)/2+\eta_1}} \dfrac{r_{2,j}\eta_2 \cos(r_{2,j}x_2)+\sin(r_{2,j}x_2)}{\sqrt{(r_{2,j}^2\eta_2^2+1)/2+\eta_2}},
  \end{align*}
  where $r_{1,i}$ is the $i$-th positive root of the characteristic equation
  $$
  (r_{1,i}^2\eta_1^2-1)\sin(r_{1,i}) = 2\eta_1r_{1,i}\cos(r_{1,i})
  $$
  in the $x_1$ direction. And the same for $r_{2,j}$, which is in the $j$-th positive root of the characteristic equation in the $x_2$ direction. After arranging the roots $r_{1,i}$ and $r_{2,j}$ in ascending order, we can immediately obtain a monotonically decreasing series of $\lambda_k$, which flavours our selection of dominant eigenfunctions in the truncated KL expansion with $R$ terms.
\end{enumerate}
For clarification, the positive indices $i$ and $j$ are mapped to index $k$ so that $\lambda_k$ is monotonically decreasing. It is noted these eigenvalues and eigenfunctions are just the product of solutions of the one-dimensional process in the Fredholm integral equation, since the covariance functions are separable.

\subsection{Training conditions}

In the training, there are different sets of parameters \rr{that} could affect the performance of the neural network. We list down as \rr{follows}:
\begin{itemize}
  \item Number of truncated terms in KL expansion: $R=4$
  \item Number of hidden layers: $L=1$
  \item Number of neurons in the hidden layer: 10
  \item Stopping criteria
  \begin{itemize}
    \item Minimum \rr{value of the} cost function gradient: $10^{-6}$ or
    \item Maximum number of training epochs: $1,000,000$
  \end{itemize}
  \item Sample size of training set: $M=10,000$
  \item Sample size of testing set: $M=500$
\end{itemize}

After we obtain an accurate neural network from the proposed algorithm, besides the NRMSE of testing set samples, we also present some characteristics of the preconditioner based on the approximate eigenvectors from the learning adaptive BDDC algorithm, which include the number of iterations required, minimum and maximum eigenvalues of the preconditioner. To show the performance of preconditioning on our predicted dominant eigenvectors, we will use the infinity norm to measure the largest difference and also the following symmetric mean absolute percentage error (sMAPE), proposed in \cite{makridakis1993accuracy}, to show a relative error type measure with an intuitive range between 0\% and 100\%:
\begin{align*}
  sMAPE = \dfrac{1}{M}\sum_{i=1}^M  \dfrac{|q_i-\widehat{q}_i|}{|q_i|+|\widehat{q}_i|},
\end{align*}
where $q_i$ are some target quantities from preconditioning and $\widehat{q}_i$ is the same type of quantity as $q_i$ from preconditioning using our \rr{predicted} eigenvectors. In the following subsections, all the computation and results were obtained using MATLAB R2019a with Intel Xeon Gold 6130 CPU and GeForce GTX 1080 Ti GPU in parallel.

\subsection{Brownian sheet covariance function}

\begin{figure}[h!]
\begin{subfigure}{.5\textwidth}
  \centering
  \includegraphics[width=1\linewidth]{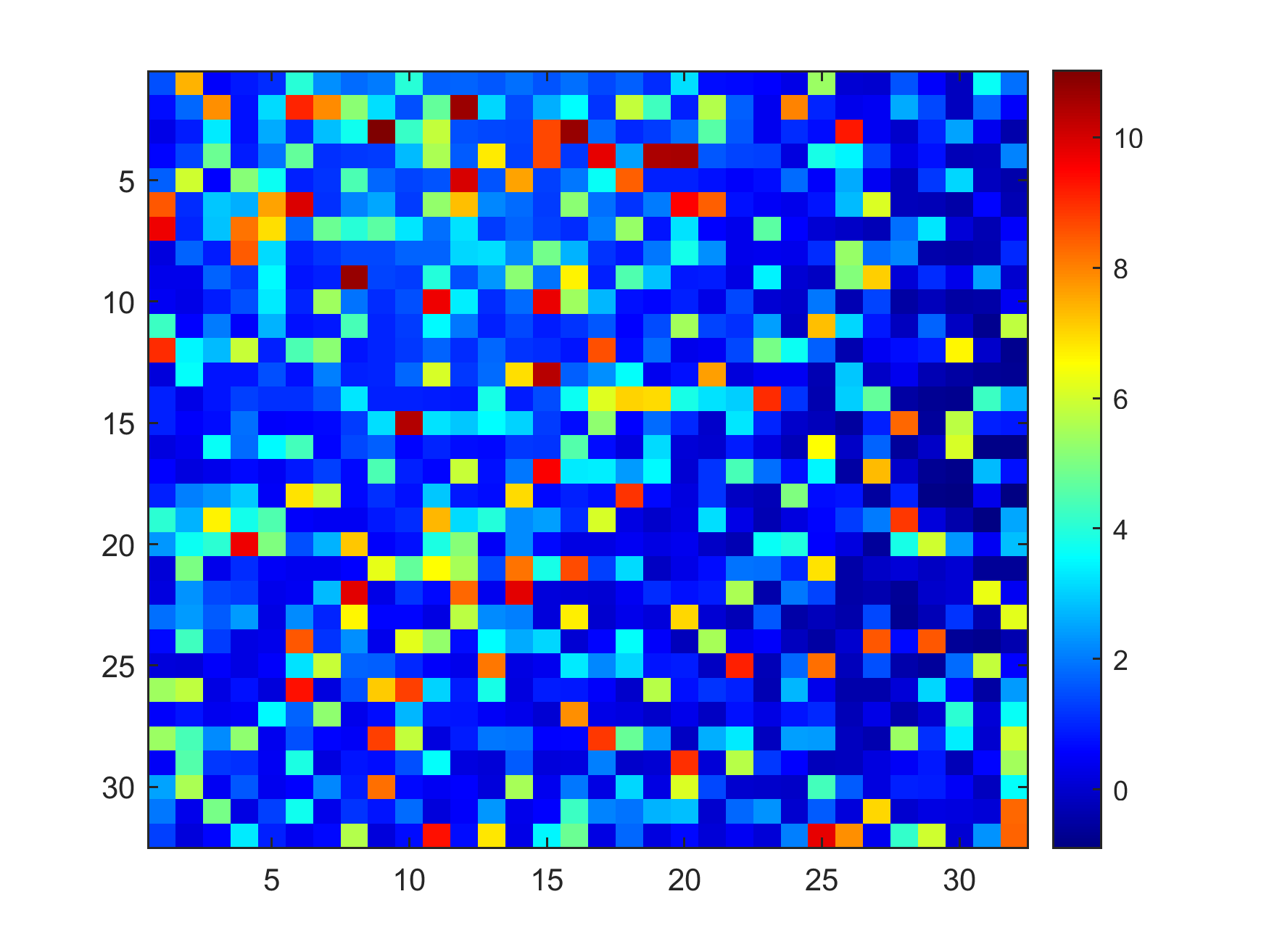}
  \caption{Expected function $\mathcal{A}$}
\end{subfigure}%
\begin{subfigure}{.5\textwidth}
  \centering
  \includegraphics[width=1\linewidth]{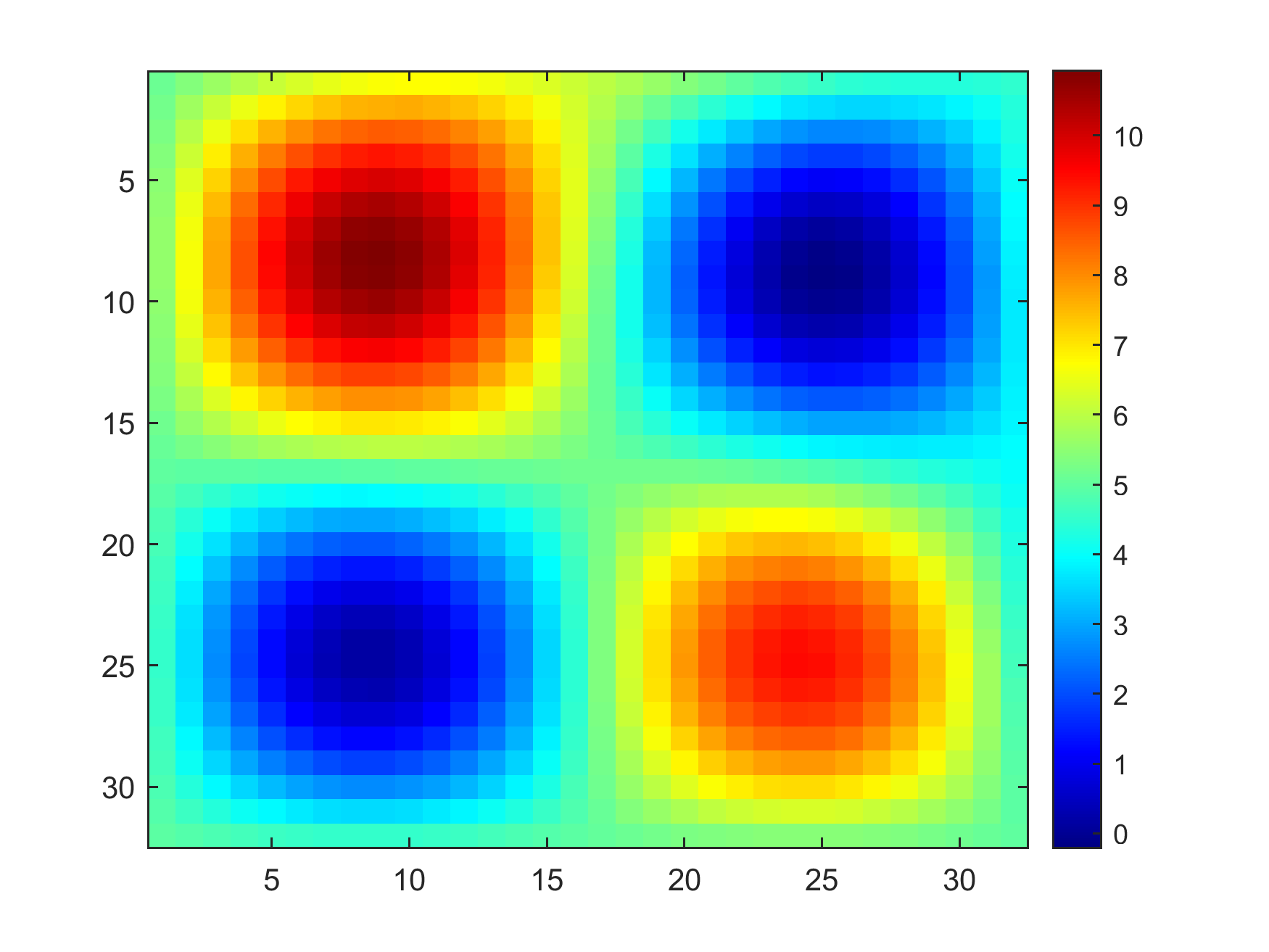}
  \caption{Expected function $\mathcal{B}$}
\end{subfigure}

\begin{subfigure}{.5\textwidth}
  \centering
  \includegraphics[width=1\linewidth]{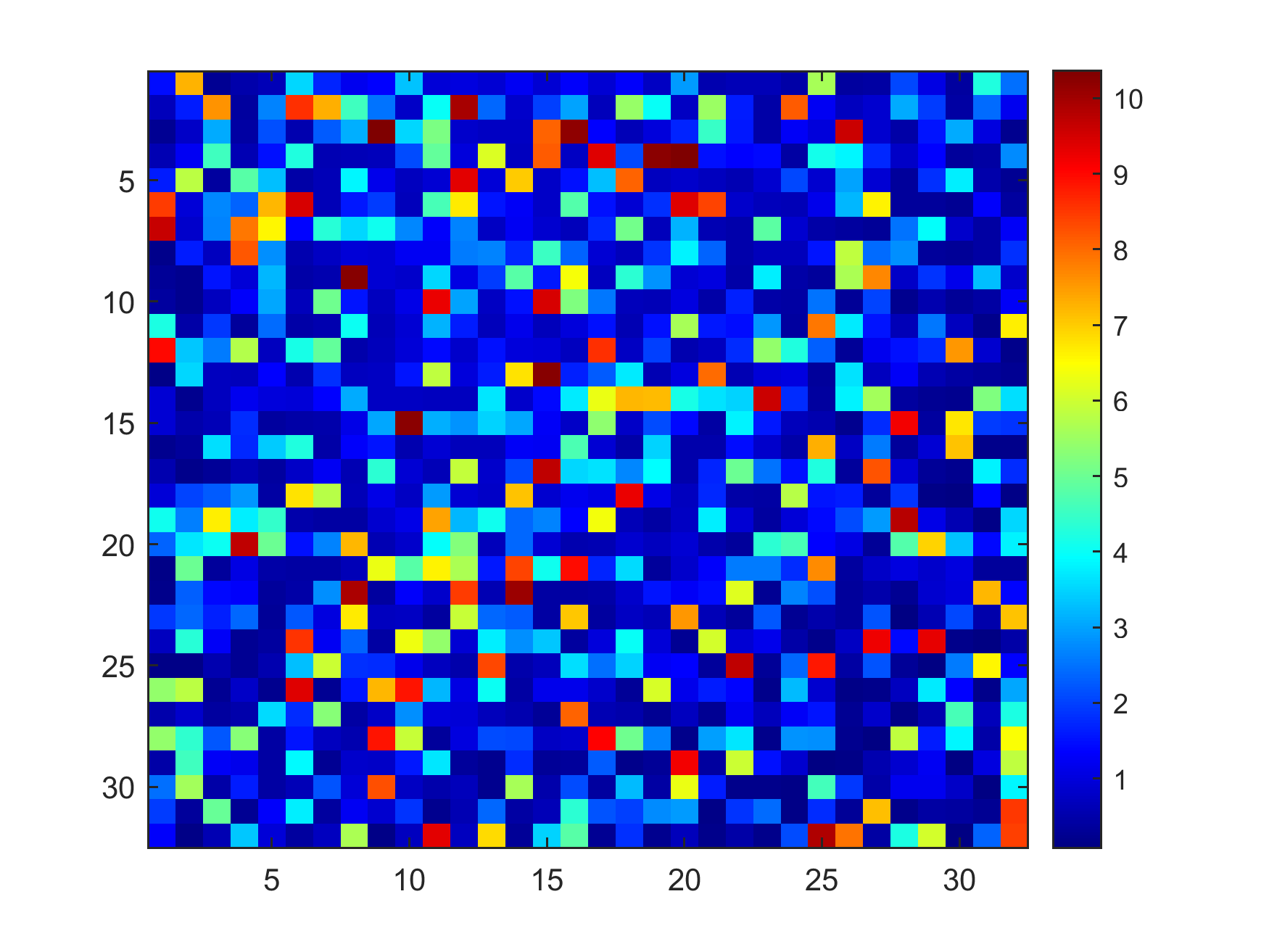}
  \caption{A realization of $K(\boldsymbol{x},\omega)$ with expected function $\mathcal{A}$}
\end{subfigure}%
\begin{subfigure}{.5\textwidth}
  \centering
  \includegraphics[width=1\linewidth]{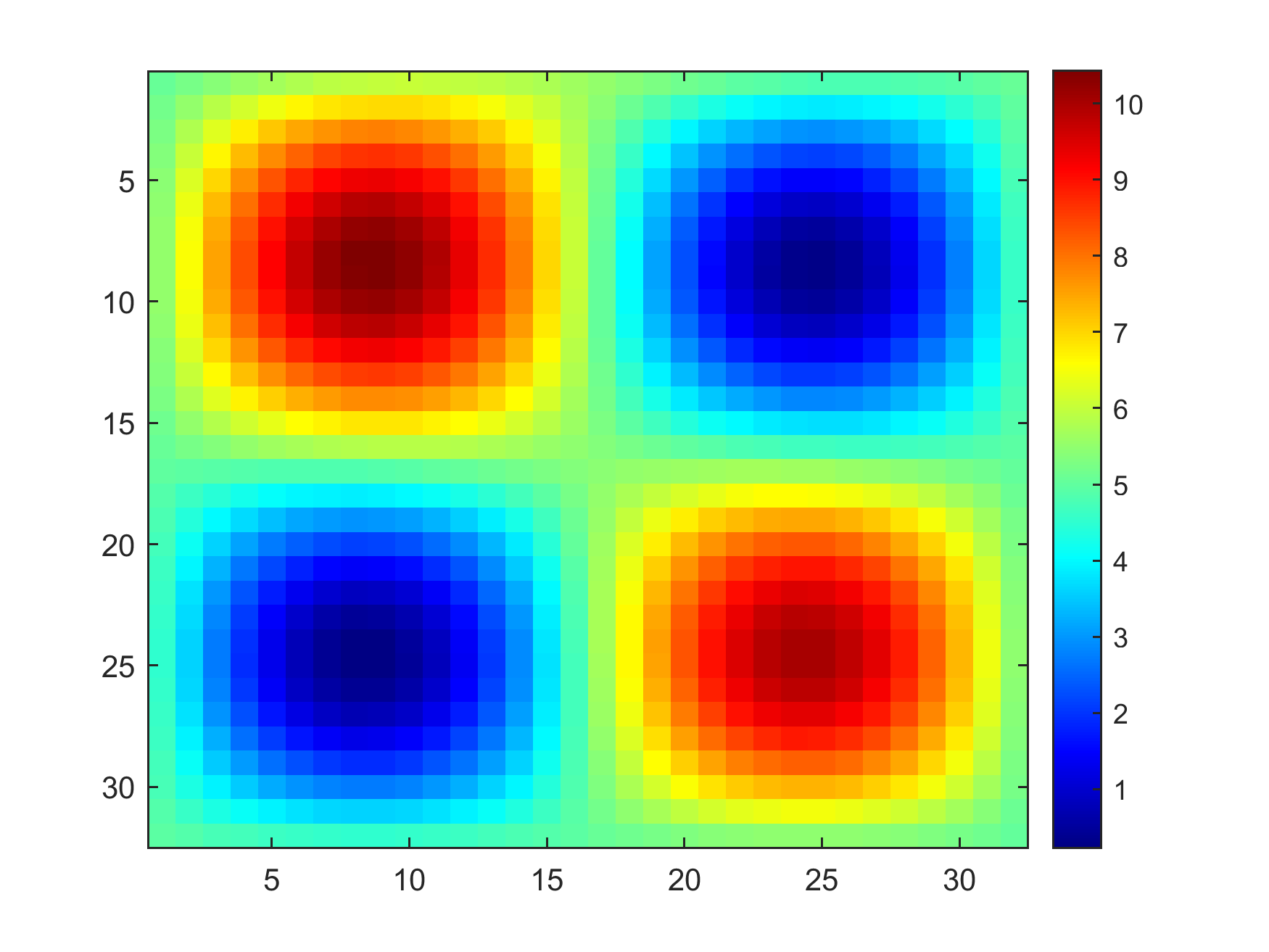}
  \caption{A realization of $K(\boldsymbol{x},\omega)$ with expected function $\mathcal{B}$}
\end{subfigure}

\begin{subfigure}{.5\textwidth}
  \centering
  \includegraphics[width=1\linewidth]{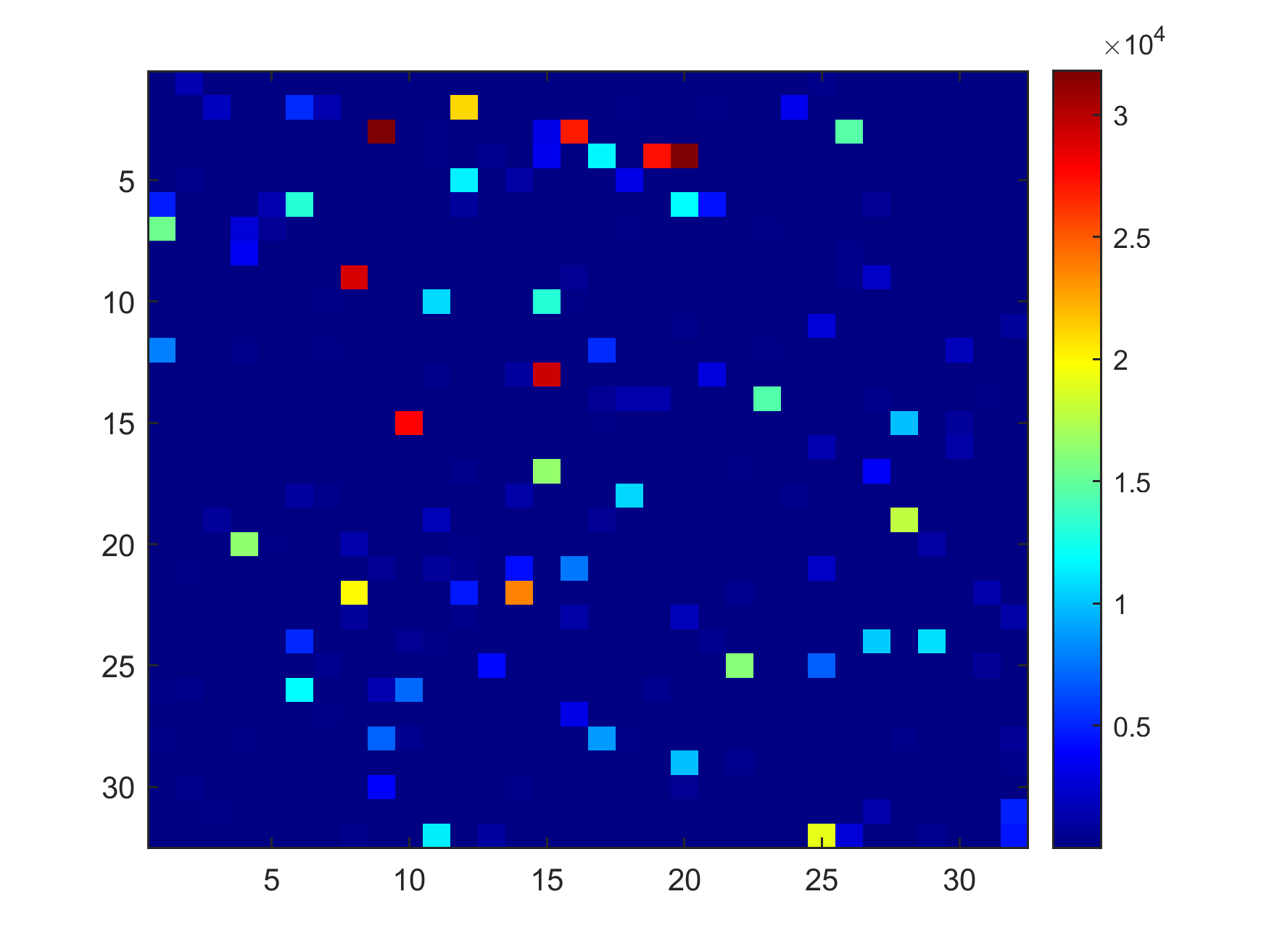}
  \caption{A realization of $\rho(\boldsymbol{x},\omega)$ with expected function $\mathcal{A}$}
\end{subfigure}%
\begin{subfigure}{.5\textwidth}
  \centering
  \includegraphics[width=1\linewidth]{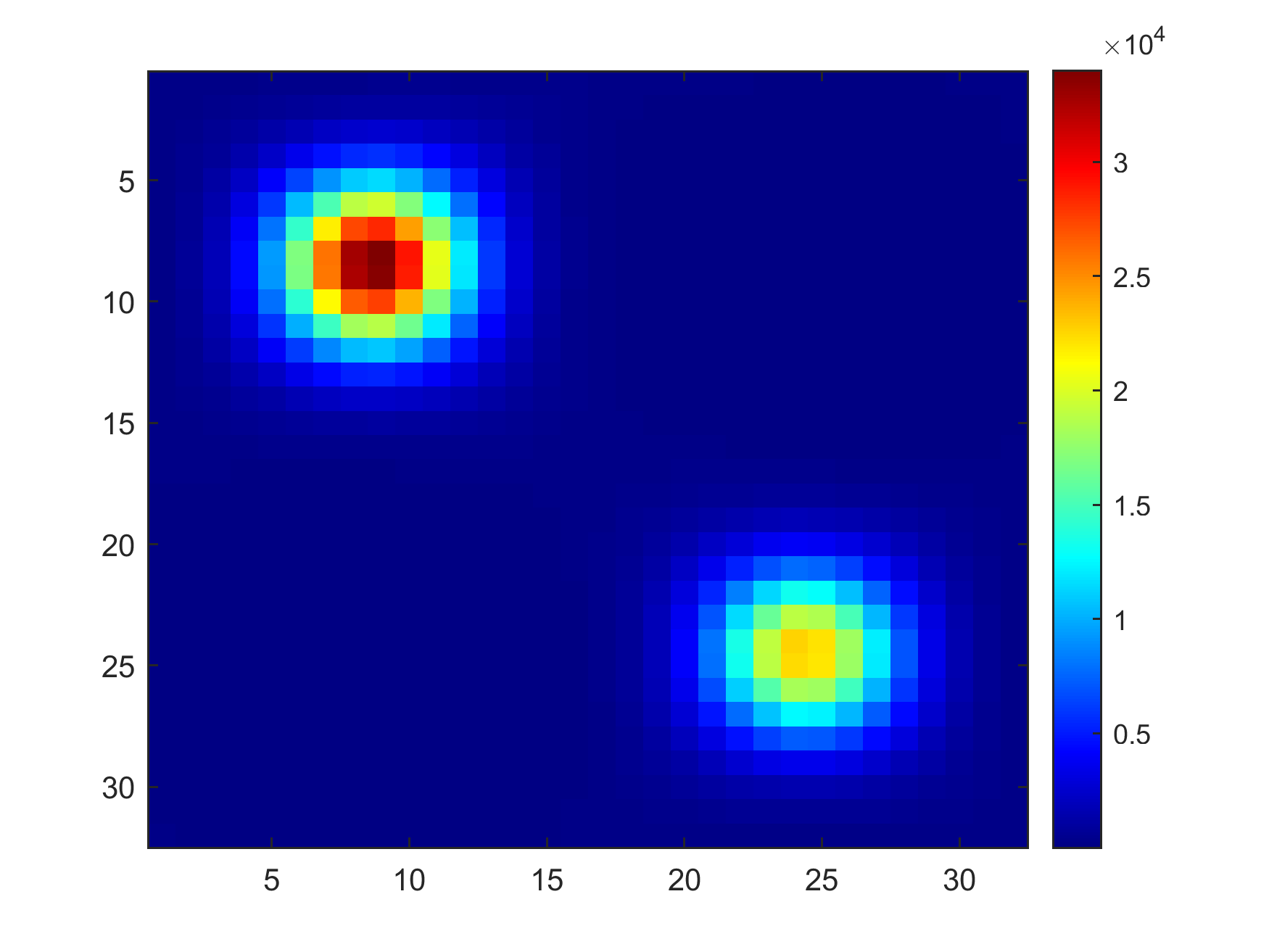}
  \caption{A realization of $\rho(\boldsymbol{x},\omega)$ with expected function $\mathcal{B}$}
\end{subfigure}
\caption{Realizations of permeability coefficient when expected functions $\mathcal{A}$ and $\mathcal{B}$ are used}
\label{fig:AB}
\end{figure}

The first set of numerical tests bases on the KL expansion of Brownian sheet covariance function with the following expected functions $E[K](\boldsymbol{x})$:
\begin{enumerate}
  \item[($\mathcal{A}$)] $10^s$
  \item[($\mathcal{B}$)] $5\sin(2\pi x_1)\sin(2\pi x_2)+5$.
\end{enumerate}
Expected function $\mathcal{A}$ is a random coefficient as $s\in[-1,1]$ is randomly chosen for each fine grid element, and expected function $\mathcal{B}$ is a smooth trigonometric function with values between $[0,10]$. The common property between these two choices is that the resulting permeability function $\rho(\boldsymbol{x},\omega)$ is highly oscillatory with magnitude order about $10^{-2}$ to $10^5$, which are good candidates to test our method on oscillatory and high contrast coefficients. Here, we show the appearance of mean functions $\mathcal{A}$ and $\mathcal{B}$ and the corresponding permeability function $\mathcal{K}(\boldsymbol{x},\omega)$ in Figure \ref{fig:AB}. We can observe that from the first row of Figure \ref{fig:AB} to the second row, that is, from the appearances of expected function $E[K](\boldsymbol{x})$ to the logarithmic permeability coefficient $K(\boldsymbol{x},\omega)$, there are no significant differences except only a minor change on the values on each fine grid element. Nevertheless, these small stochastic changes cause a high contrast $\rho(\boldsymbol{x},\omega)$ after taking exponential function as in the third row of Figure \ref{fig:AB}.

In Table \ref{train:AB}, the training results of these two expected functions $\mathcal{A}$ and $\mathcal{B}$ under the training conditions mentioned are presented. Although the random initial choice of conjugate gradient \rr{algorithm} for minimization may have influences, in general, the neural network of
expected function $\mathcal{A}$ usually needs more epochs until stopping criterion is reached. However, this increase in training epochs does not bring a better training NRMSE when compared with expected function $\mathcal{B}$. We can see the same phenomenon when the testing set is considered.
The main reason is due to the smoothness of expected function $\mathcal{B}$, the neural network can capture its characteristics easier than a random coefficient.

\begin{table}[h]
\centering
\begin{tabular}{|c|c|c|}\hline
 \rr{Case}  & Epochs trained & Training NRMSE \\ \hline
 $\mathcal{A}$   &   2.29e+05  &  1.97e-02    \\
 $\mathcal{B}$   &   5.39e+04  &  3.77e-03    \\ \hline
\end{tabular}
\caption{Training records \rr{when} expected functions $\mathcal{A}$ and $\mathcal{B}$ are considered}
\label{train:AB}
\end{table}

In the testing results, we consider another set of data containing $M=500$ samples as a testing set. After obtaining the predicted eigenvectors from neural network, we plug \rr{it} into the BDDC \rr{preconditioned} solver and obtain an estimated preconditioner. \rr{For the estimated preconditioner, we report} several properties such as the number of iteration needed for \rr{the iterative solver}, the maximum and minimum eigenvalues of the \rr{preconditioned system}. Therefore, to better show the performance of network, besides the testing error NRMSE, we also use the sMAPE and $l_\infty$ norm of difference of the iteration numbers, the maximum and minimum eigenvalues of the estimated preconditioner and the target preconditioner. The comparison results are listed below:

\begin{table}[h]
\centering
\begin{tabular}{|c|c|c|c|c|c|c|c|}\hline
 \multirow{2}{*}{\rr{Case}} & \multirow{2}{*}{Testing NRMSE} & \multicolumn{3}{c|} {sMAPE ( $l_\infty$ error) in}  \\ \cline{3-5}
  & &  Iteration  number & $\lambda_{min}$ & $\lambda_{max}$   \\ \hline
 $\mathcal{A}$  & 6.58e-02 & 6.94e-02 ( 2 ) &  2.58e-07 ( 4.02e-06 ) &  4.45e-03 ( 3.58e-02 )    \\
 $\mathcal{B}$  & 6.73e-03 & 6.93e-02 ( 1 ) &  3.19e-05 ( 1.50e-04 ) &  2.01e-03 ( 4.94e-02 )     \\  \hline
\end{tabular}
\caption{Testing records \rr{when} expected functions $\mathcal{A}$ and $\mathcal{B}$ are considered}
\label{test:AB}
\end{table}

It is clear that the testing NRMSE of expected function $\mathcal{B}$ is also much smaller than expected function $\mathcal{A}$. One possible reason is because the function used is smooth, and the neural network can better capture the property of the resulting stochastic permeability function. We can see that values of the testing NRMSE of both expected function $\mathcal{A}$ and expected function $\mathcal{B}$ are just a bit larger than the training NRMSE. This suggests \rr{that} the neural network considered is capable to give \rr{a} good prediction on dominant eigenvectors in the coarse space even when the magnitude of coefficient function value changes dramatically across each fine grid element. Nevertheless, even though a more accurate set of eigenvectors can be obtained for \rr{the} expected function $\mathcal{B}$, but \rr{its} performance after the BDDC preconditioning is not always better than that of expected function $\mathcal{A}$. Obviously, $l_\infty$ errors in the aspects of $\lambda_{max}$ and $\lambda_{min}$ are larger \rr{for the expected function $\mathcal{B}$}.

\begin{figure}[h!]
\vspace{1cm}
\begin{subfigure}{.5\textwidth}
  \centering
  \includegraphics[width=1\linewidth]{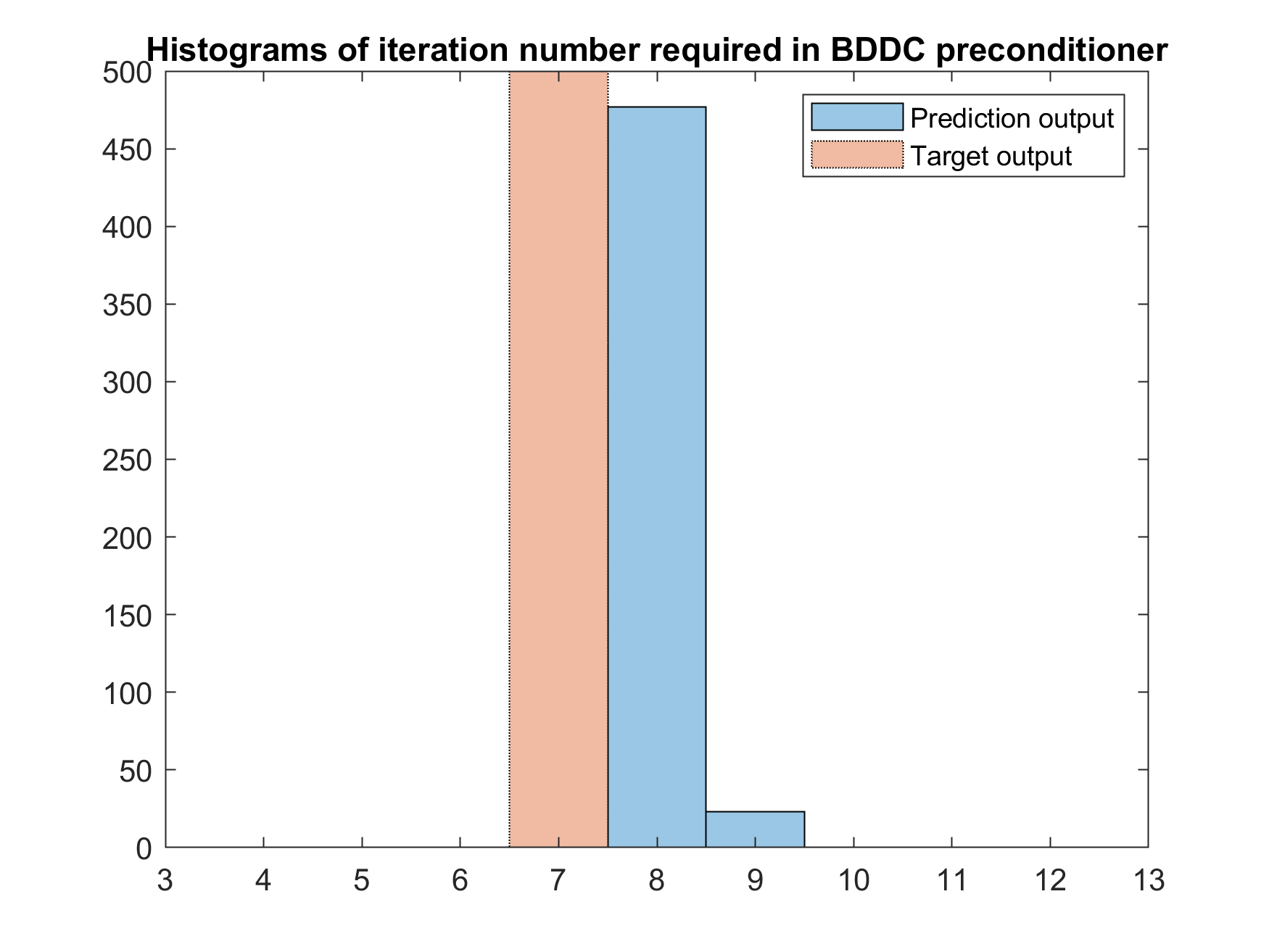}
\end{subfigure}%
\begin{subfigure}{.5\textwidth}
  \centering
  \includegraphics[width=1\linewidth]{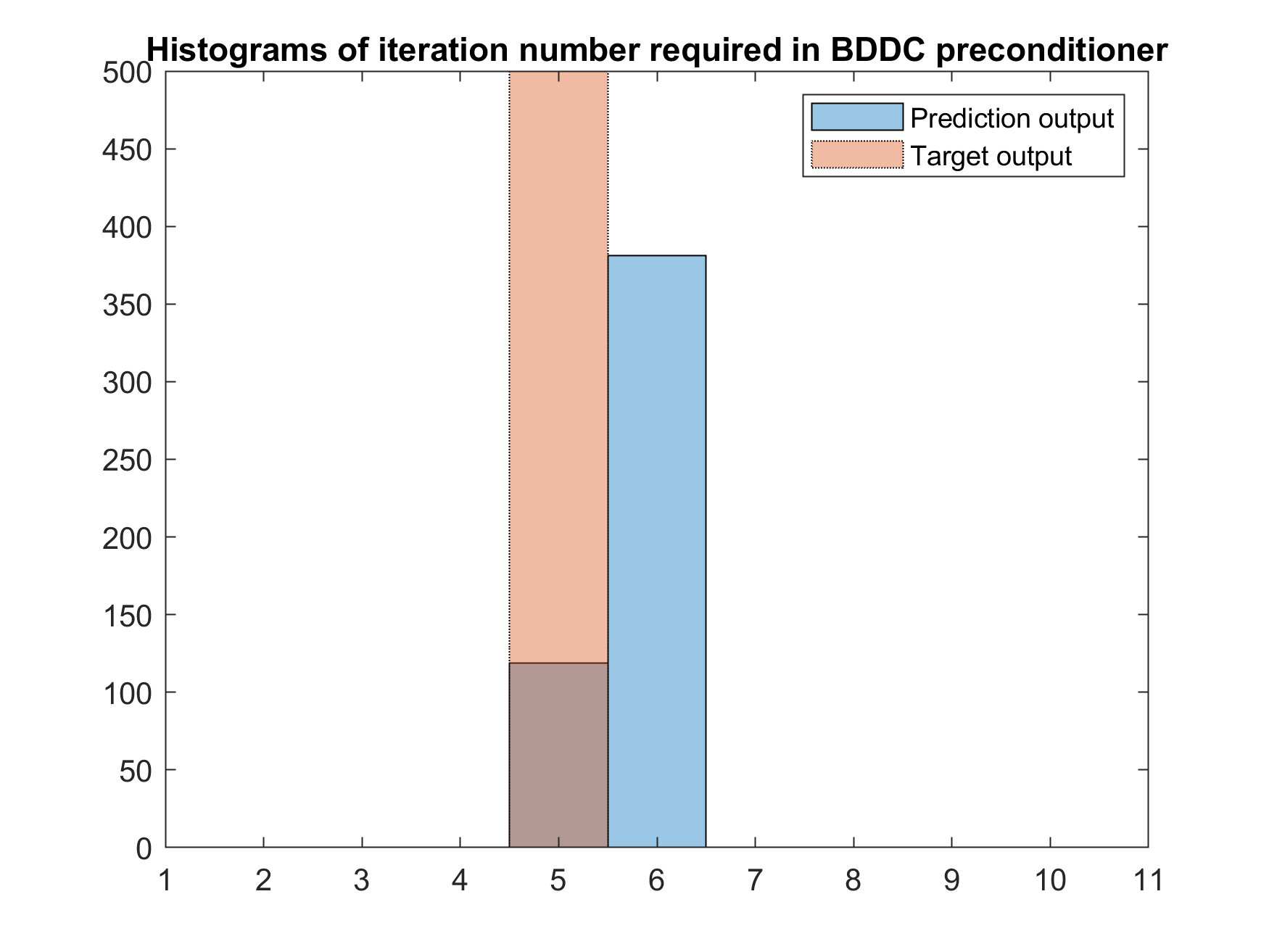}
\end{subfigure}

\begin{subfigure}{.5\textwidth}
  \centering
  \includegraphics[width=1\linewidth]{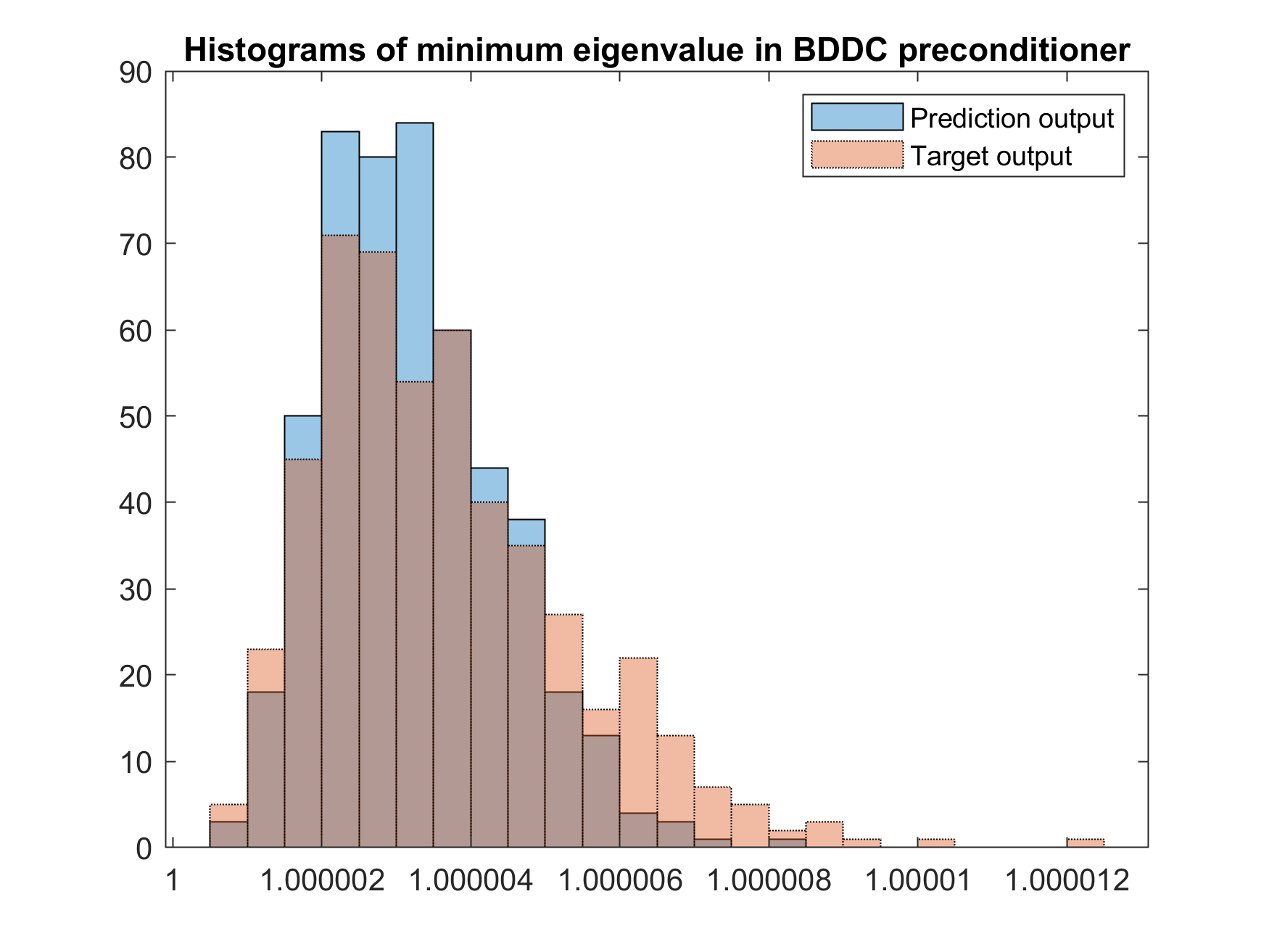}
\end{subfigure}%
\begin{subfigure}{.5\textwidth}
  \centering
  \includegraphics[width=1\linewidth]{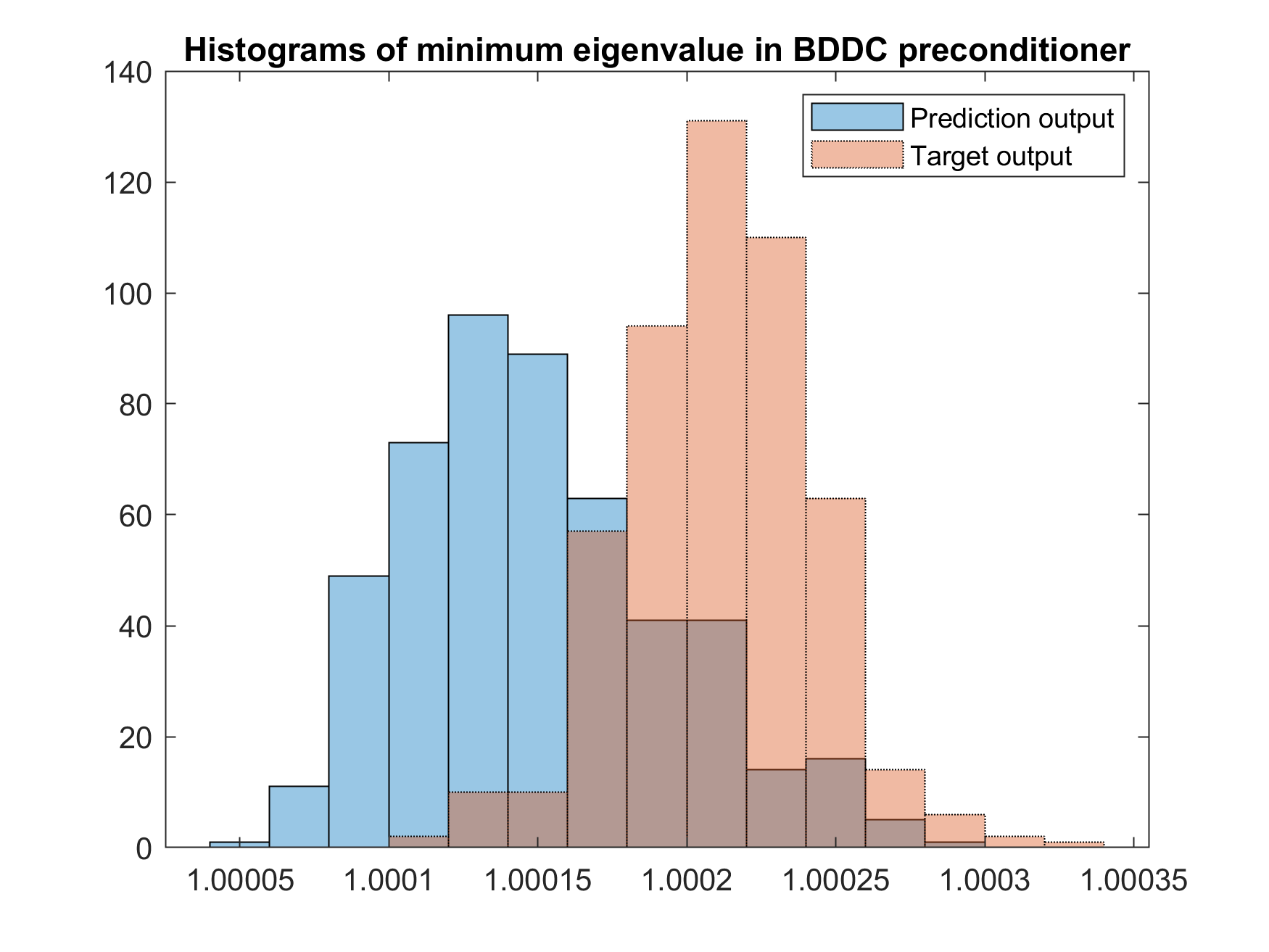}
\end{subfigure}

\begin{subfigure}{.5\textwidth}
  \centering
  \includegraphics[width=1\linewidth]{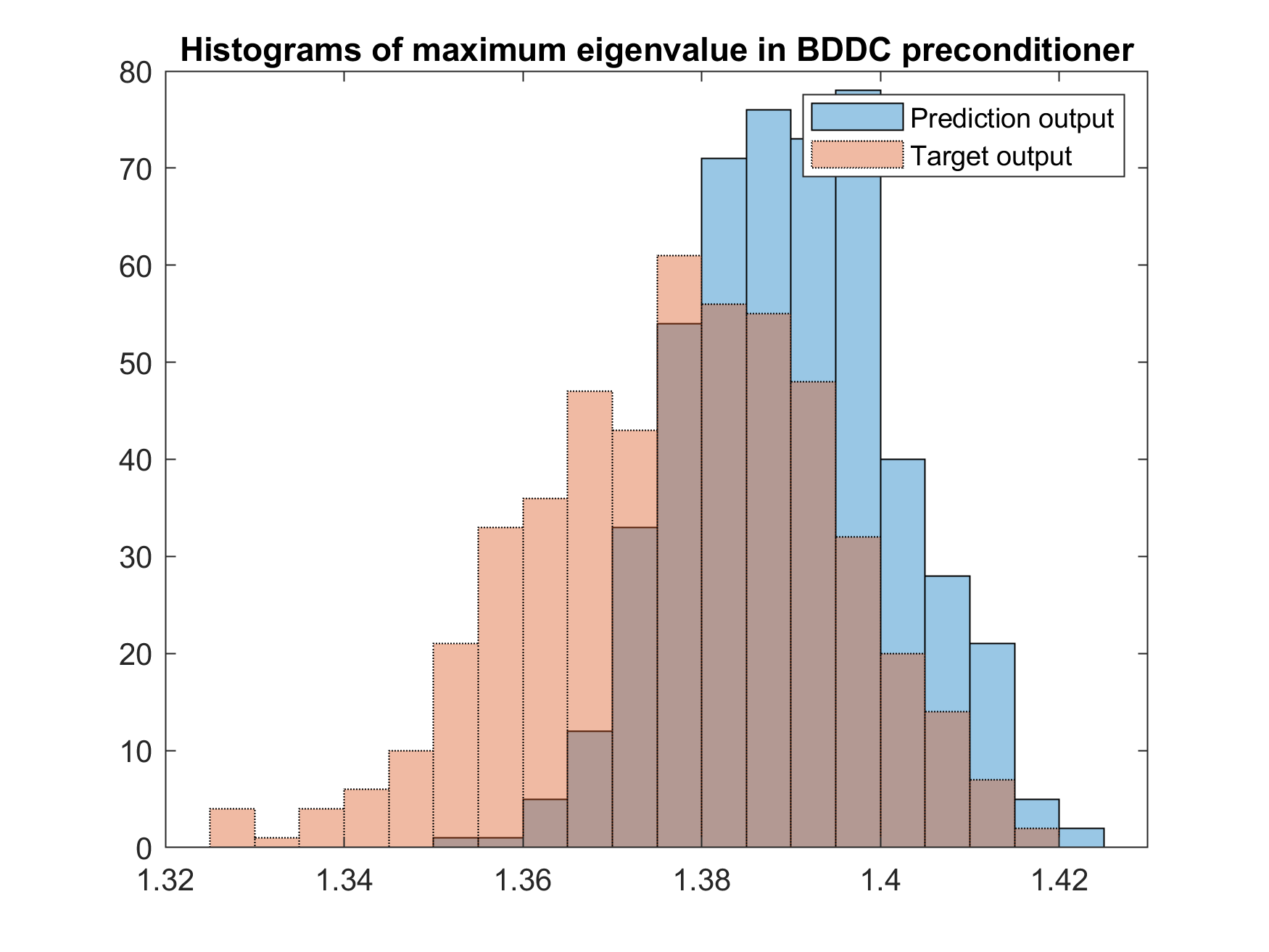}
\end{subfigure}%
\begin{subfigure}{.5\textwidth}
  \centering
  \includegraphics[width=1\linewidth]{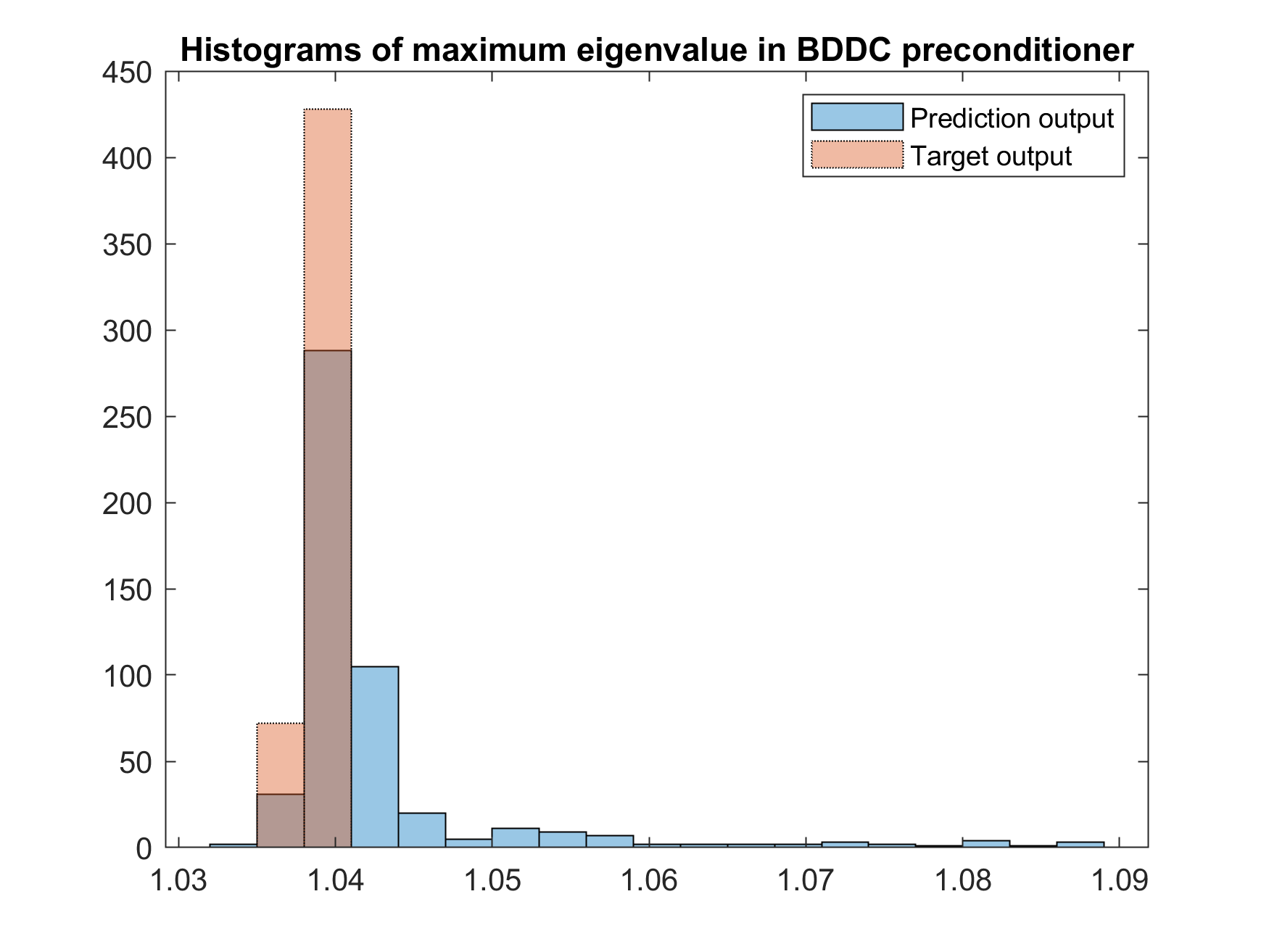}
\end{subfigure}
\caption{Comparisons in \rr{the performance of the} BDDC preconditioner when expected functions $\mathcal{A}$ (left column) and $\mathcal{B}$ (right column) are used}
\label{fig:ABcompare}
\end{figure}

Detailed comparisons can be found in Figure \ref{fig:ABcompare}, \rr{where} the area of overlapping represents how close is the specified quantity of the estimated preconditioner compared to the target preconditioner. In the first row of figures, we can see the number of iteration required for \rr{the predicted} eigenvectors is usually more than the target iteration number, however, the differences in the iteration number is not highly related to the difference in minimum and maximum eigenvalues. We can see more examples to verify it in this and later results. The bin width of the histogram of  minimum eigenvalues for expected function $\mathcal{A}$ is about $10^{-7}$, which is in coherence with the $l_\infty$ error, but we eliminate this small magnitude effect and the sMAPE is considered. The high ratio of overlapping area for expected function $\mathcal{A}$ confirms the usage of sMAPE as a good measure, as the resulting sMAPE 2.57e-07 of expected function $\mathcal{A}$ is much smaller than 3.19e-05 of expected function $\mathcal{B}$.

For the last column of errors, the $l_\infty$ norm of maximum eigenvalues is larger for expected function $\mathcal{B}$. Nevertheless, it does not mean \rr{that the} neural network of expected function $\mathcal{B}$ performs worse. From its histogram, we can see an extremely concentrated and overlapping area at the bin around 1.04, moreover, there are some outliers for the prediction results, which is one of the source for a larger $l_\infty$ norm. Therefore, the sMAPE can well represent the performance of preconditioning on the estimated results, and we can conclude that both neural networks show a good results on oscillatory and high contrast coefficients, which can be represented by NRMSE and sMAPE. In the next subsection, besides an artificial coefficient, we will consider one that is closer to \rr{a} real life case.

Although the $l_\infty$ norm seems \rr{not describe} the characteristics of performance in a full picture, it is still a good and intuitive measure, for example, in the difference of iteration number, we can immediately realize how large will \rr{be} the worst case.


%
%

\subsection{Exponential covariance function}

\begin{figure}[h!]
\begin{subfigure}{.5\textwidth}
  \centering
  \includegraphics[width=1\linewidth]{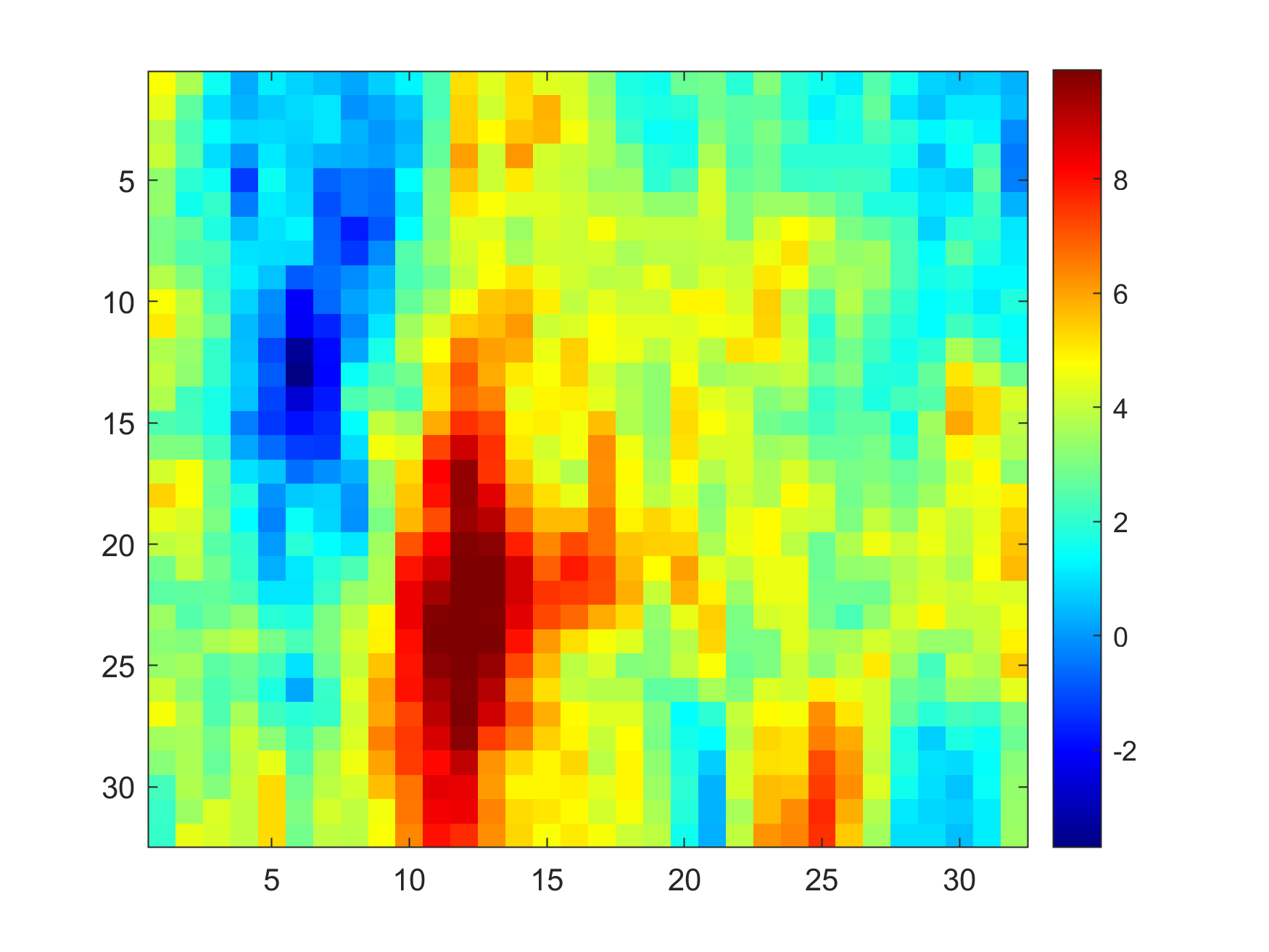}
  \caption{Mean permeability field of Layer 34}
\end{subfigure}
\begin{subfigure}{.5\textwidth}
  \centering
  \includegraphics[width=1\linewidth]{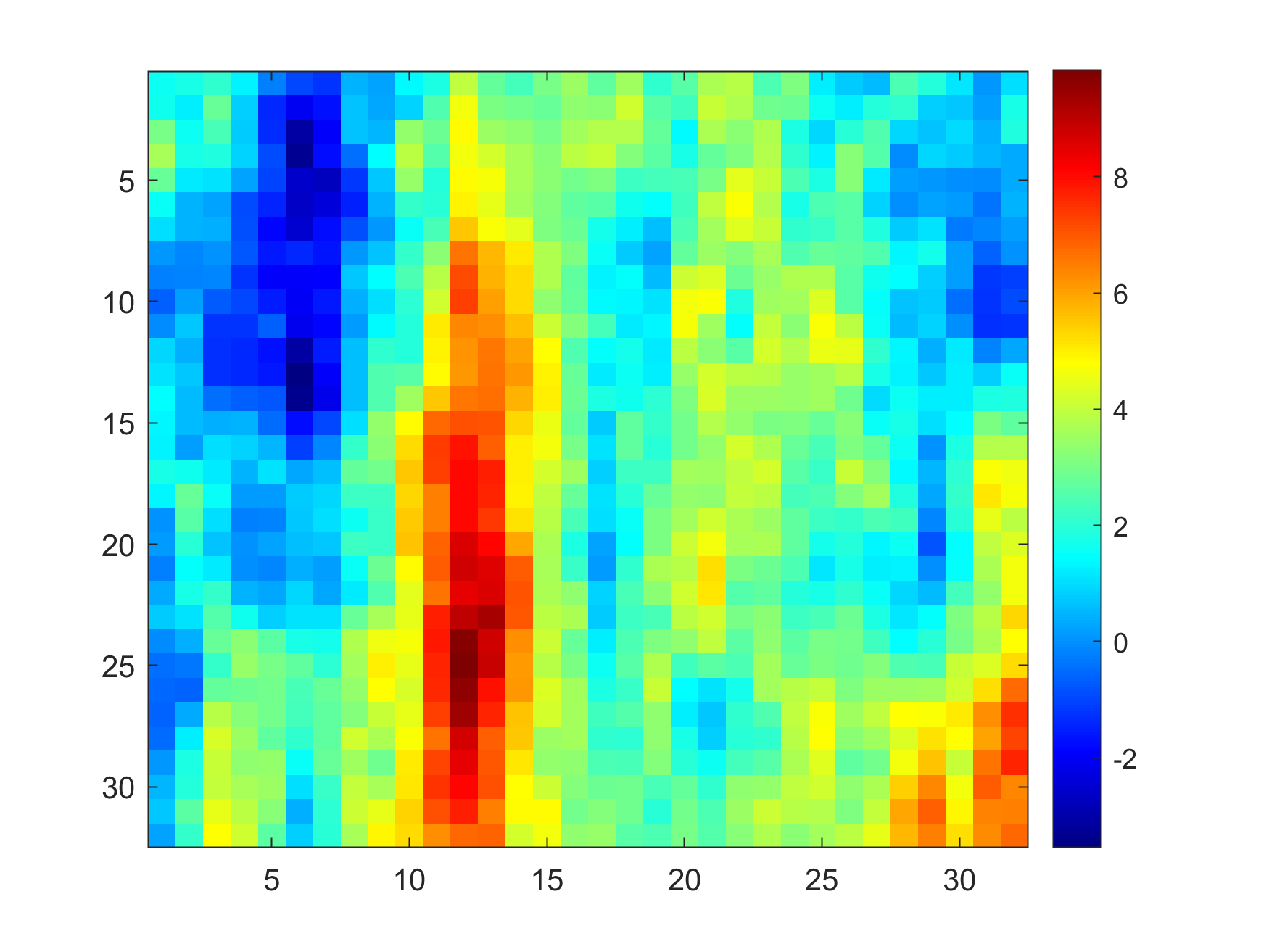}
  \caption{Mean permeability field of Layer 35}
\end{subfigure}%
\caption{Mean permeability fields of different layers in the $x$-direction}
\label{fig:diffmean}
\end{figure}

To test our method on realistic highly varying coefficients, we consider expected functions that come from the second model of the 10th SPE Comparative Solution Project (SPE10), with the KL expanded exponential covariance function as the stochastic source. For clarification, in this set of experiment, we first use the modified permeability of Layer 35 in the $x$-direction of SPE10 data as the expected function $E[K](\boldsymbol{x})$ with $(\sigma_K^2,\eta_1,\eta_2)=(1,0.25,0.25)$ in the exponential covariance function $C_K(\boldsymbol{x},\hat{\boldsymbol{x}})$. We \rr{then} train the neural network  and the resulting neural network is used to test its generalization capacity on the following different situations:
\begin{enumerate}
  \item Different stochastic behaviour:
        \begin{itemize}
          \item All remain unchanged except \rr{$(\eta_1,\eta_2)=(0.2,0.125)$}
        \end{itemize}
  \item Different mean permeability function:
        \begin{itemize}
          \item The expected function is changed to the Layer 34 in the $x$-direction of SPE10 data\rr{, but the stochastic parameters are unchanged, i.e., $(\sigma_K^2,\eta_1,\eta_2)=(1,0.25,0.25)$}
        \end{itemize}
\end{enumerate}
Note, all the permeability fields used are modified into our computational domain and shown in Figure \ref{fig:diffmean} for \rr{a} clear comparison. In the Layer 35 permeability field, four sharp properties can be observed. There are one blue strip on the very left side, and a reversed but narrower red strip is next to it. In the top right corner, a low permeability area is shown and just below it, a small high permeability area is located in the bottom right corner. So a similar but with slightly different properties of Layer 34 are chosen to test the generalization capacity of our neural network.

\begin{figure}[h!]
\begin{subfigure}{.5\textwidth}
  \centering
  \includegraphics[width=1\linewidth]{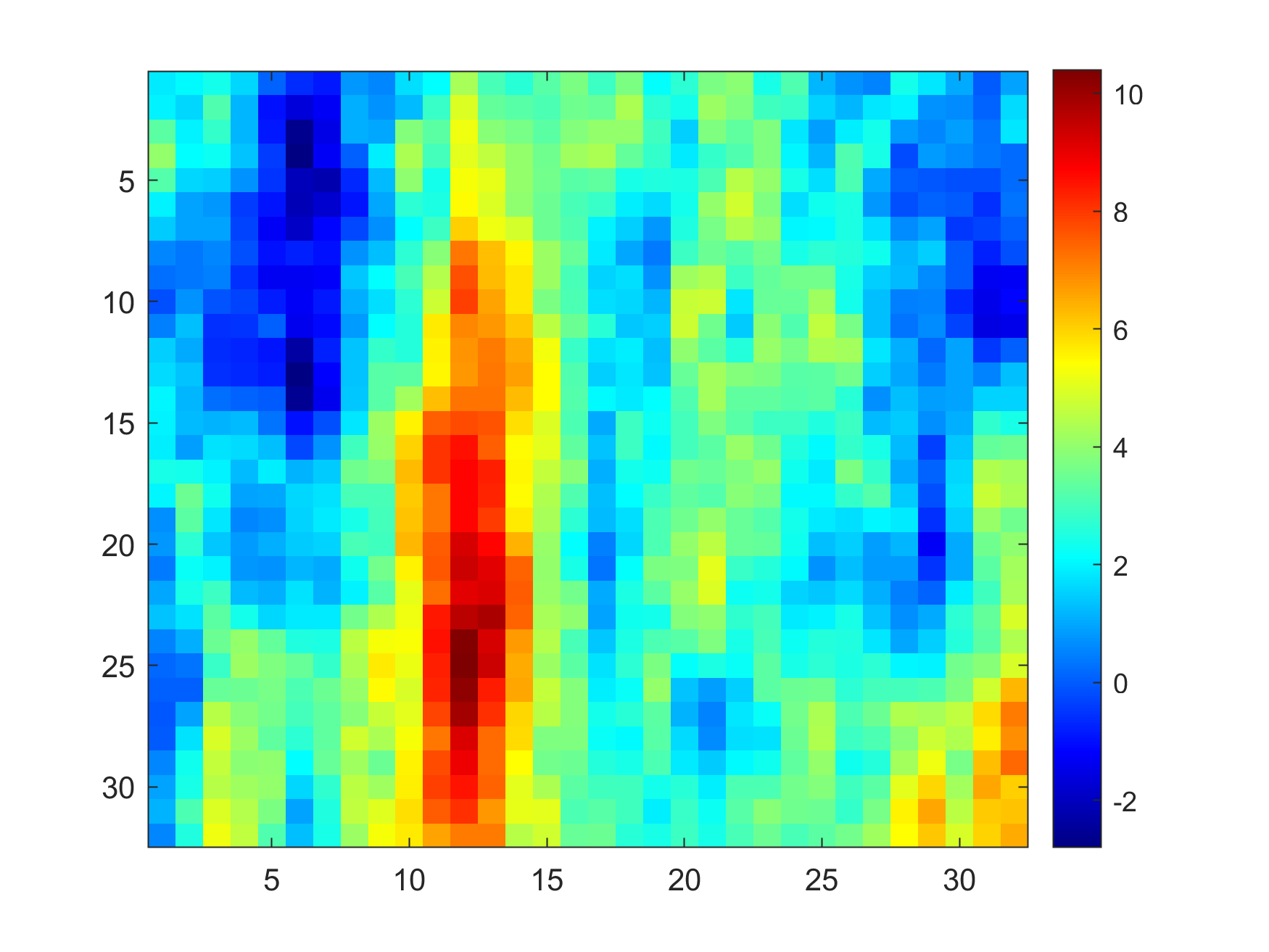}
  \caption{Log permeability when $(\sigma_K^2,\eta_1,\eta_2)=(1,0.25,0.25)$}
\end{subfigure}%
\begin{subfigure}{.5\textwidth}
  \centering
  \includegraphics[width=1\linewidth]{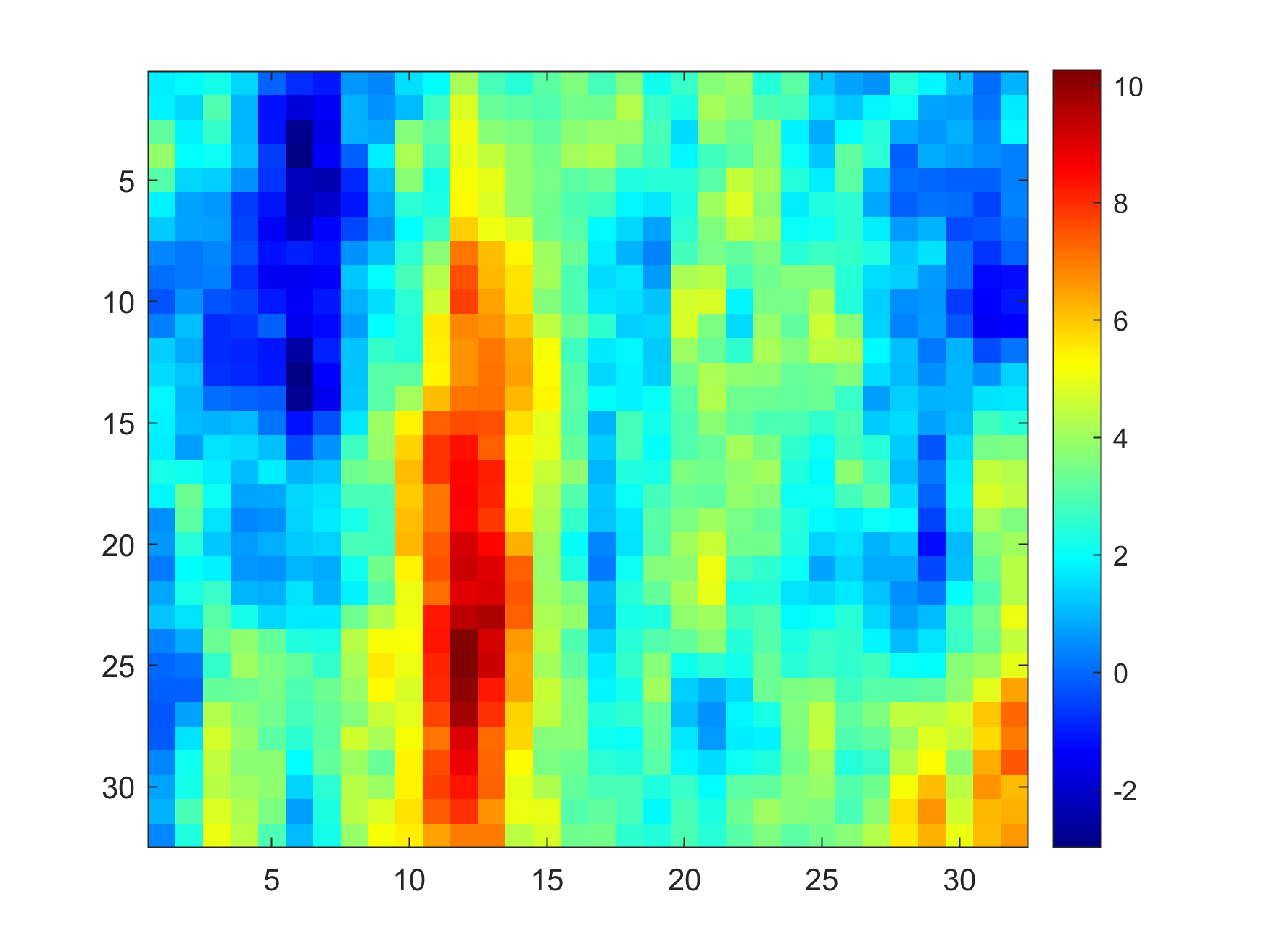}
  \caption{Log permeability  when $(\sigma_K^2,\eta_1,\eta_2)=(1,0.2,0.125)$}
\end{subfigure}
\begin{subfigure}{.5\textwidth}
  \centering
  \includegraphics[width=1\linewidth]{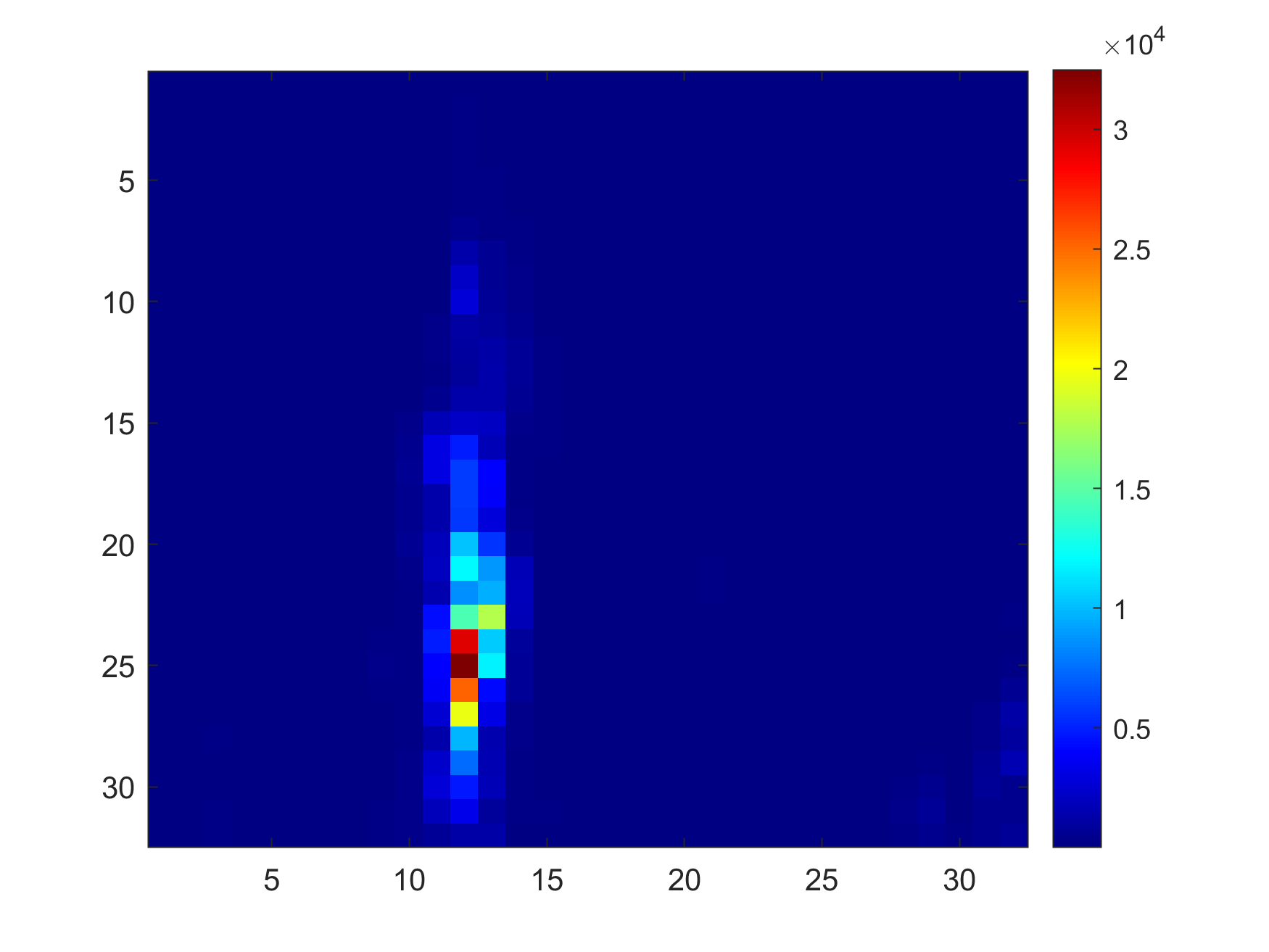}
  \caption{Permeability $\rho(\boldsymbol{x},\omega)$ when $(\sigma_K^2,\eta_1,\eta_2)=(1,0.25,0.25)$}
\end{subfigure}%
\begin{subfigure}{.5\textwidth}
  \centering
  \includegraphics[width=1\linewidth]{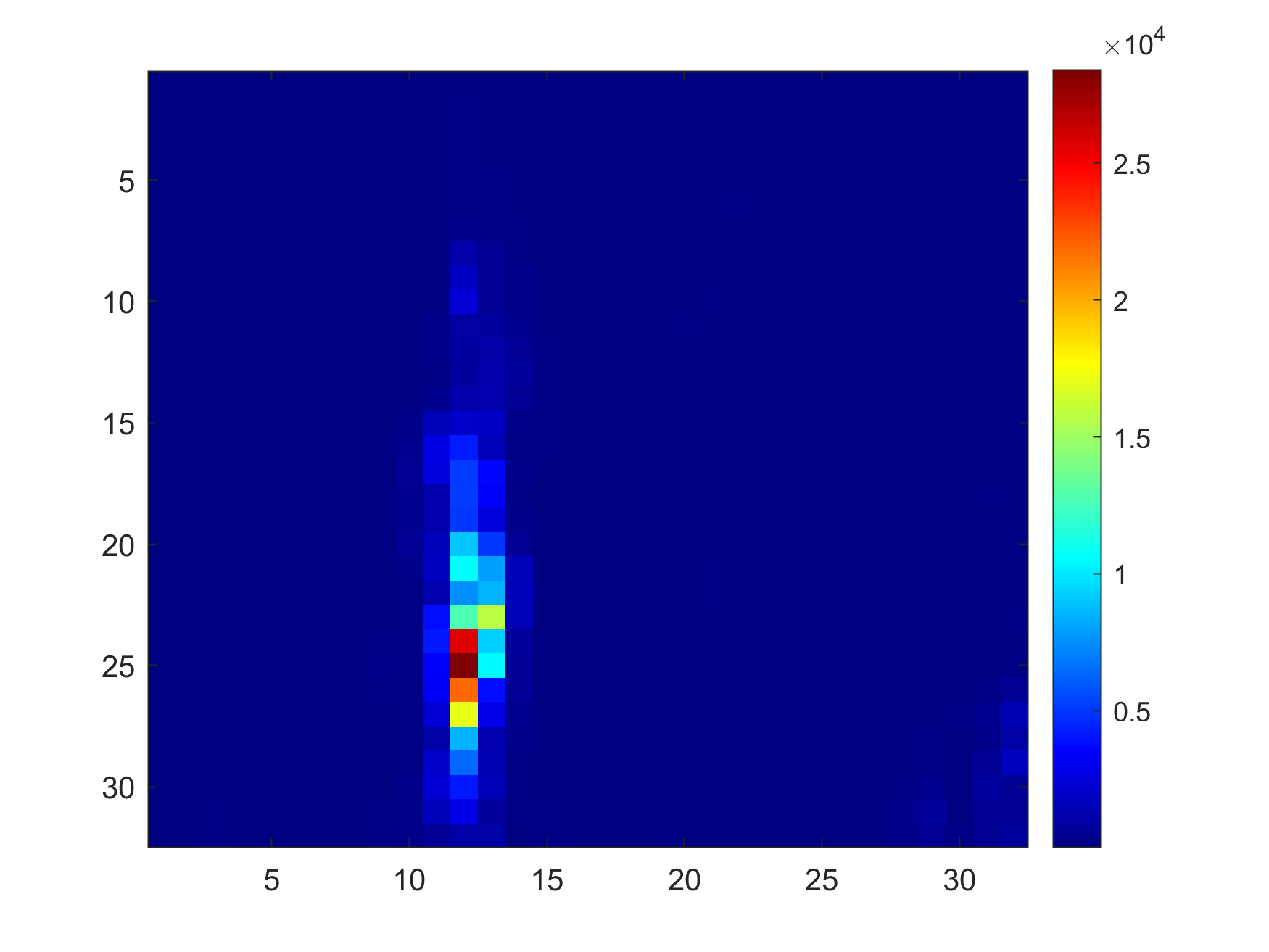}
  \caption{Permeability $\rho(\boldsymbol{x},\omega)$ when $(\sigma_K^2,\eta_1,\eta_2)=(1,0.2,0.125)$}
\end{subfigure}%
\caption{Realizations of \rr{Layer 35} permeability coefficient when different $\eta_1$ and $\eta_2$ are used}
\label{fig:diffperm}
\end{figure}

Before showing the training and testing records of neural network, we further present the logarithmic permeability $K(\boldsymbol{x},\omega)$ and permeability $\rho(\boldsymbol{x},\omega)$ when $(\sigma_K^2,\eta_1,\eta_2)=(1,0.25,0.25)$ and $(\sigma_K^2,\eta_1,\eta_2)=(1,0.2,0.125)$ in Figure \ref{fig:diffperm}. The sets of $\{\xi_i\}$ used in these two column of realizations are the same. We can \rr{thus see the} change in parameters does cause the stochastic behaviour to change, too. We present the training and testing records of neural network as below in Table \ref{train:Layer35} and Table \ref{test:Layer35}, \rr{where we specify the expected functions considered in the first column.}

\begin{table}[h]
\centering
\begin{tabular}{|c|c|c|}\hline
 \rr{Case} & Epochs trained & Training NRMSE \\ \hline
 Layer 35   &  2.62e+05  &  1.52e-02    \\ \hline
\end{tabular}
\caption{Training record when Layer 35 permeability is considered}
\label{train:Layer35}
\end{table}


\begin{table}[h!]
\centering
\setlength\tabcolsep{6.5pt}
\begin{tabular}{|c|c|c|c|c|c|c|c|}\hline
  \multirow{2}{*}{\rr{Case}} & \multirow{2}{*}{Testing NRMSE} & \multicolumn{3}{c|} {sMAPE ( $l_\infty$ error) in}  \\ \cline{3-5}
  & &  Iteration  number & $\lambda_{min}$ & $\lambda_{max}$   \\ \hline
 Layer 35  & 2.48e-02 & 3.65e-02 ( 1 ) &  7.52e-06 ( 7.50e-05 ) &  1.36e-03 ( 6.66e-02 )    \\
 \rr{Layer $35^*$}  & 2.27e-02 & 3.89e-02 ( 1 ) &  7.57e-06 ( 7.17e-05 ) &  1.27e-03 ( 6.59e-02 )    \\
 Layer 34  & 5.08e-02 & 2.52e-03 ( 1 ) &  3.25e-05 ( 1.74e-04 ) &  3.06e-02 ( 8.47e-02 )    \\ \hline
\end{tabular}
\caption{Testing records when Layer 34, 35 and \rr{$35^*$} are considered}
\label{test:Layer35}
\end{table}


\rr{In Table \ref{test:Layer35}, we clarify that the results of Layer $35^*$ also use Layer 35 permeability as the expected function but with $(\eta_1,\eta_2)=(0.2,0.125)$. Each row corresponds to the results of different testing sets, however, all the results are obtained using the same neural network of Layer 35 in Table \ref{train:Layer35}. Therefore, the testing results of Layer 35 are used as a reference to decide whether the results of Layer $35^*$ and Layer 34 are good or not.}

\begin{figure}[h!]
\vspace{1cm}
\begin{subfigure}{.5\textwidth}
  \centering
  \includegraphics[width=1\linewidth]{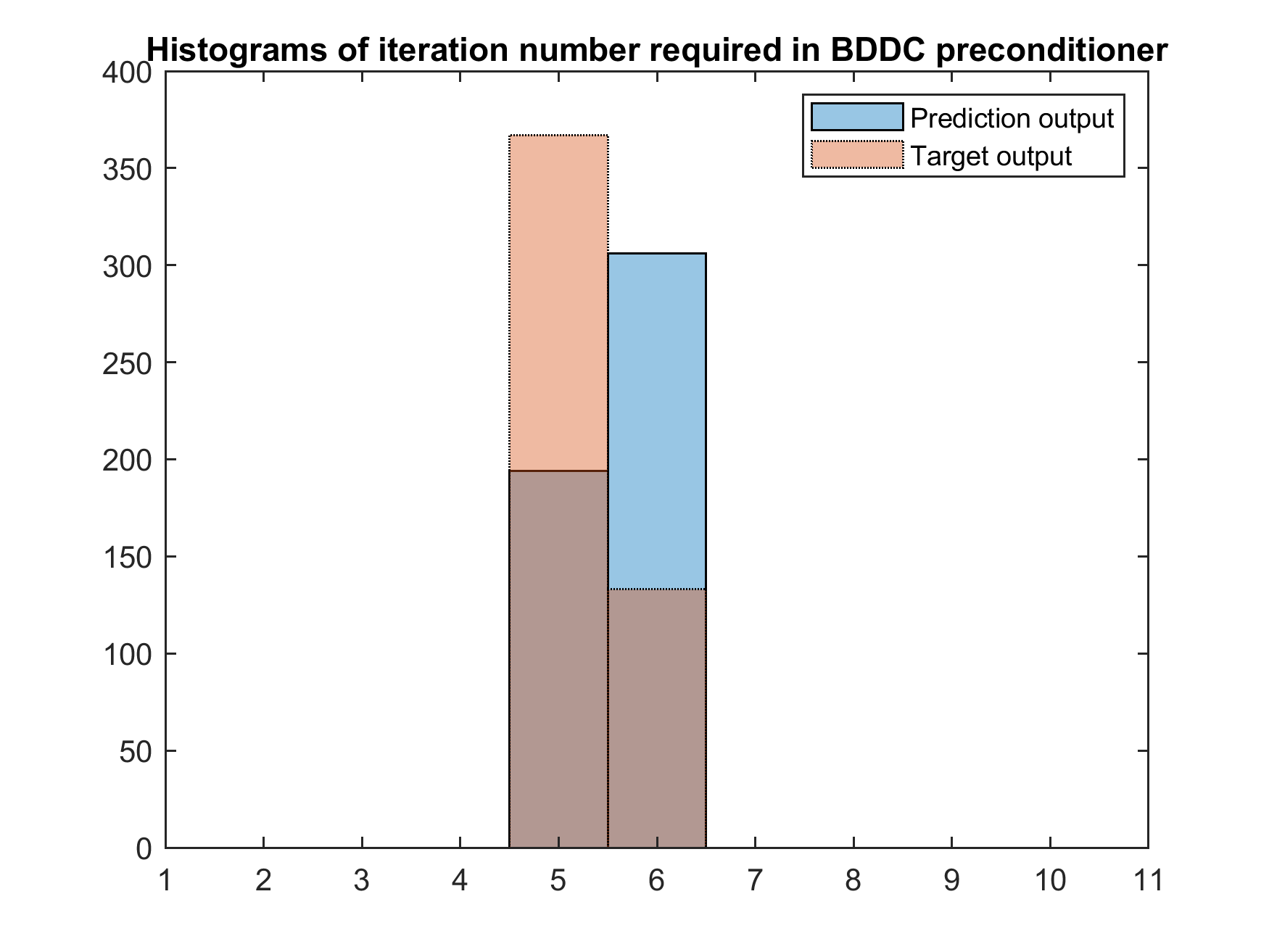}
\end{subfigure}%
\begin{subfigure}{.5\textwidth}
  \centering
  \includegraphics[width=1\linewidth]{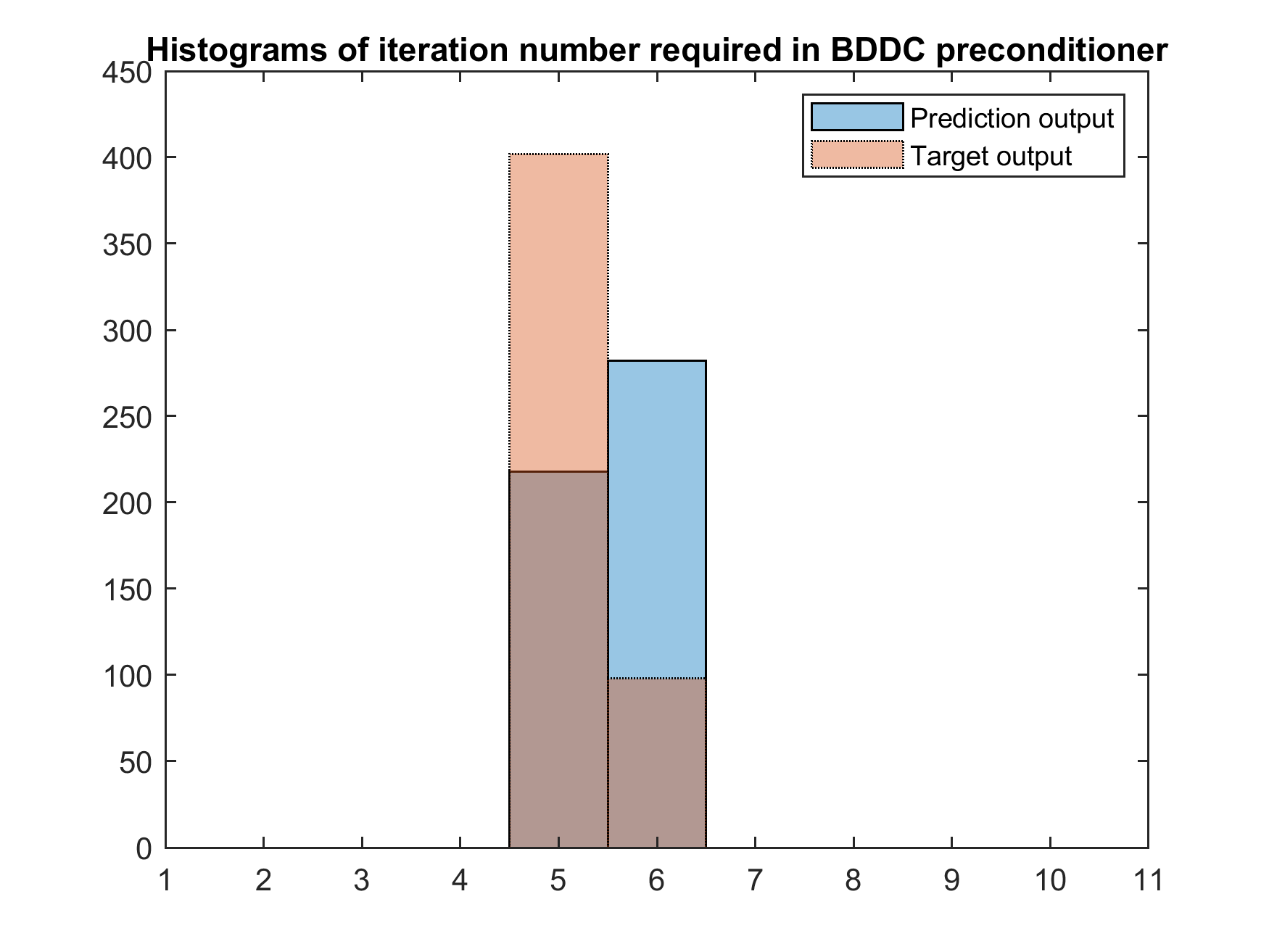}
\end{subfigure}

\begin{subfigure}{.5\textwidth}
  \centering
  \includegraphics[width=1\linewidth]{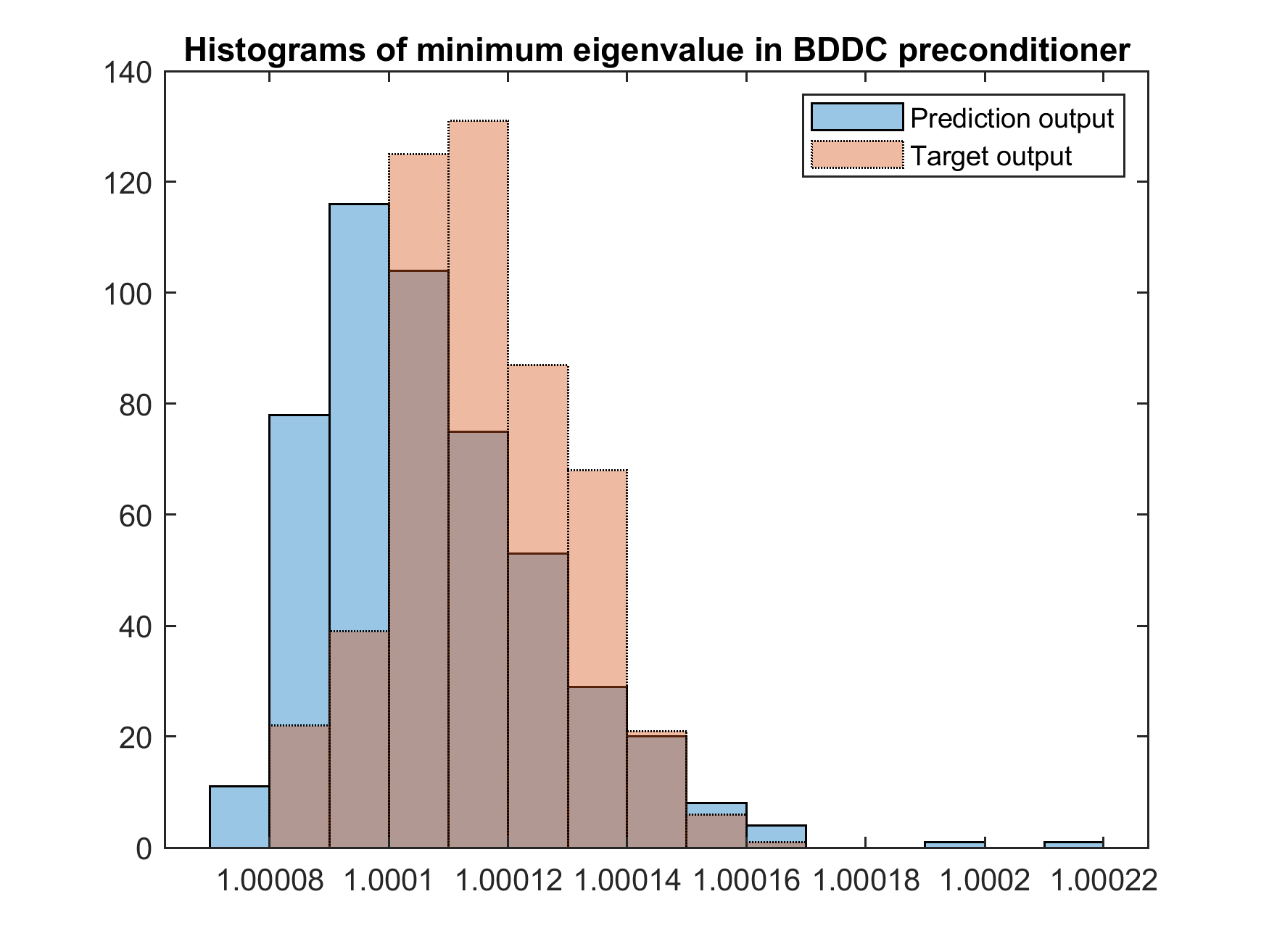}
\end{subfigure}%
\begin{subfigure}{.5\textwidth}
  \centering
  \includegraphics[width=1\linewidth]{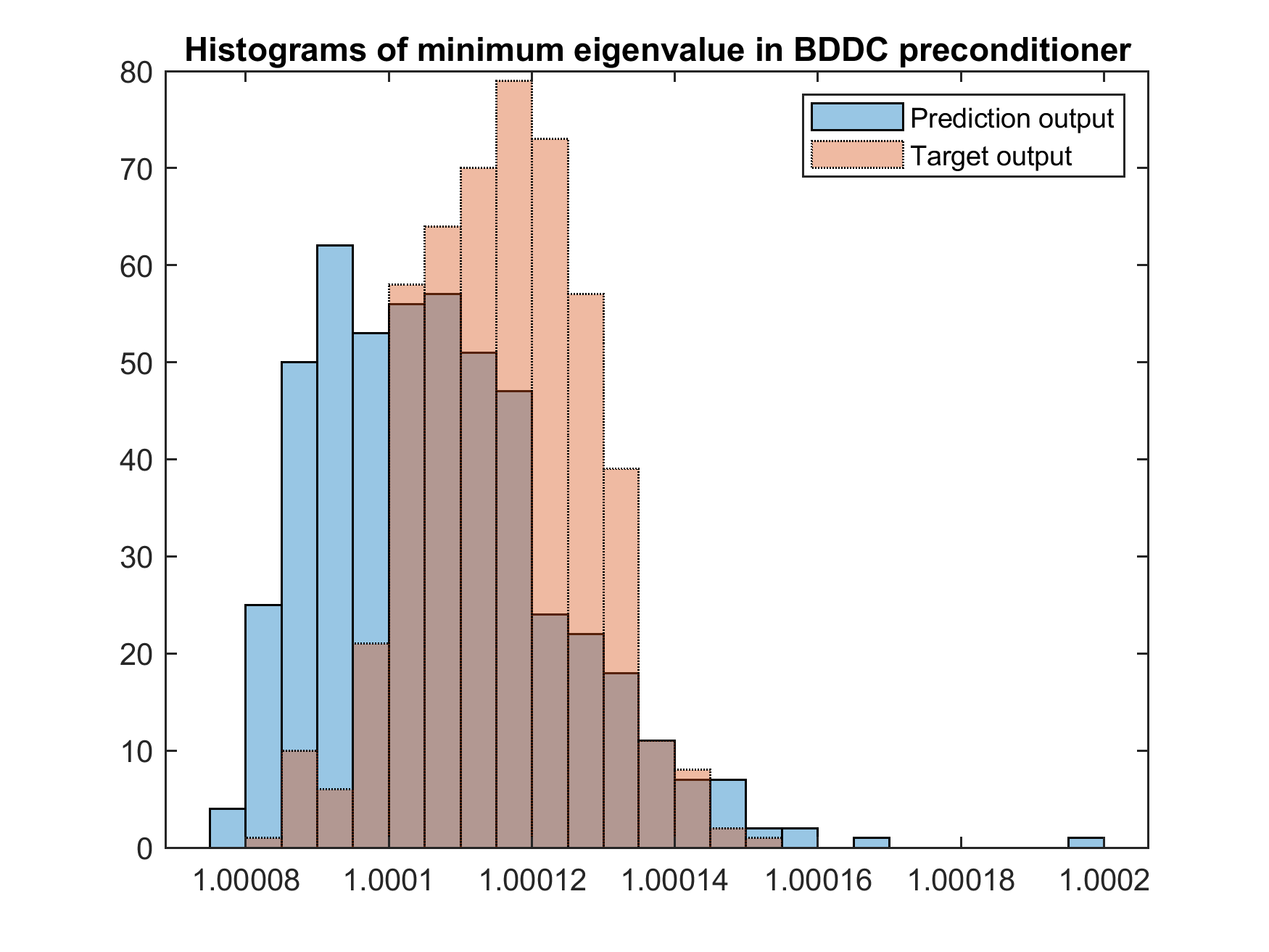}
\end{subfigure}

\begin{subfigure}{.5\textwidth}
  \centering
  \includegraphics[width=1\linewidth]{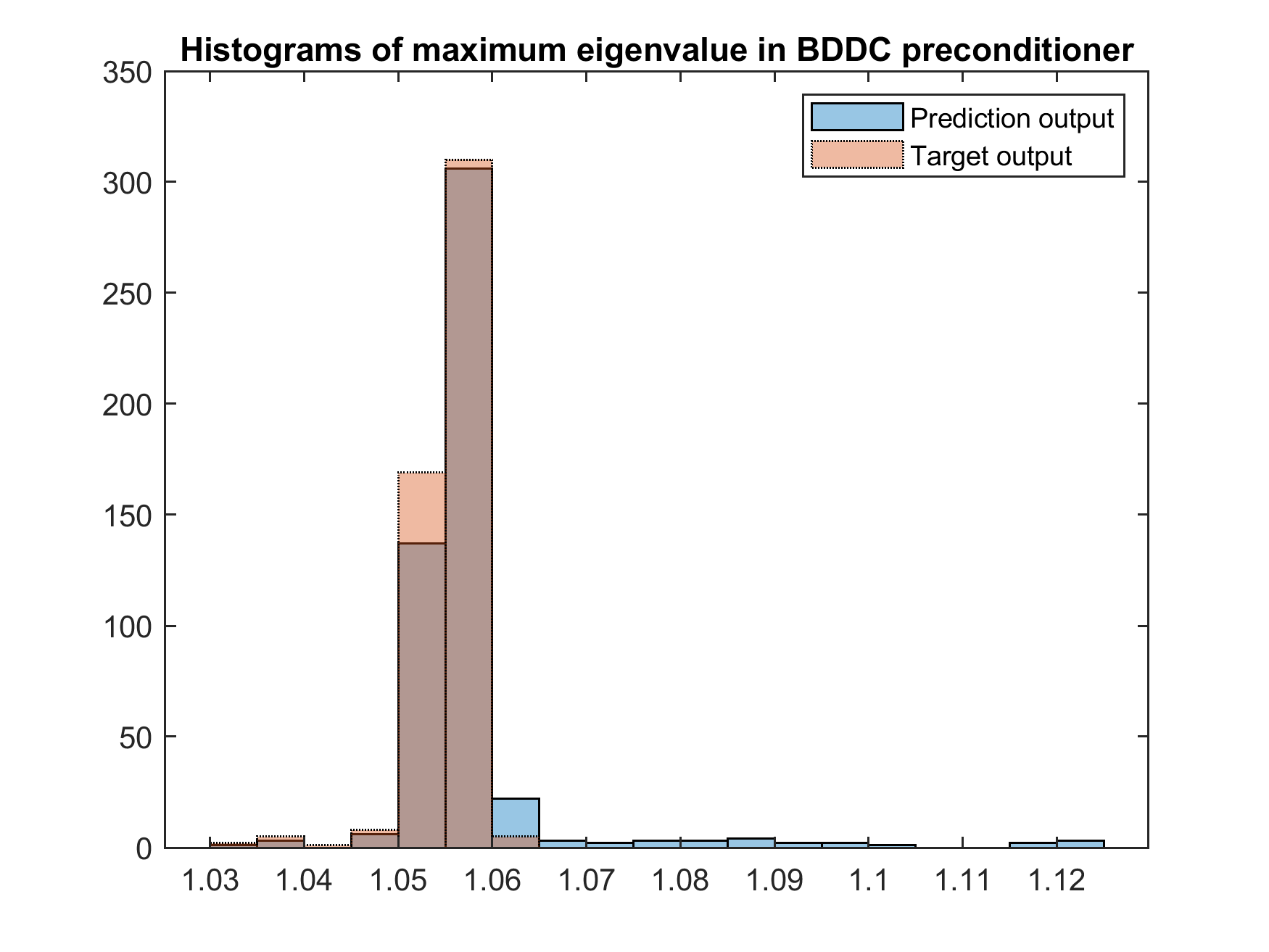}
\end{subfigure}%
\begin{subfigure}{.5\textwidth}
  \centering
  \includegraphics[width=1\linewidth]{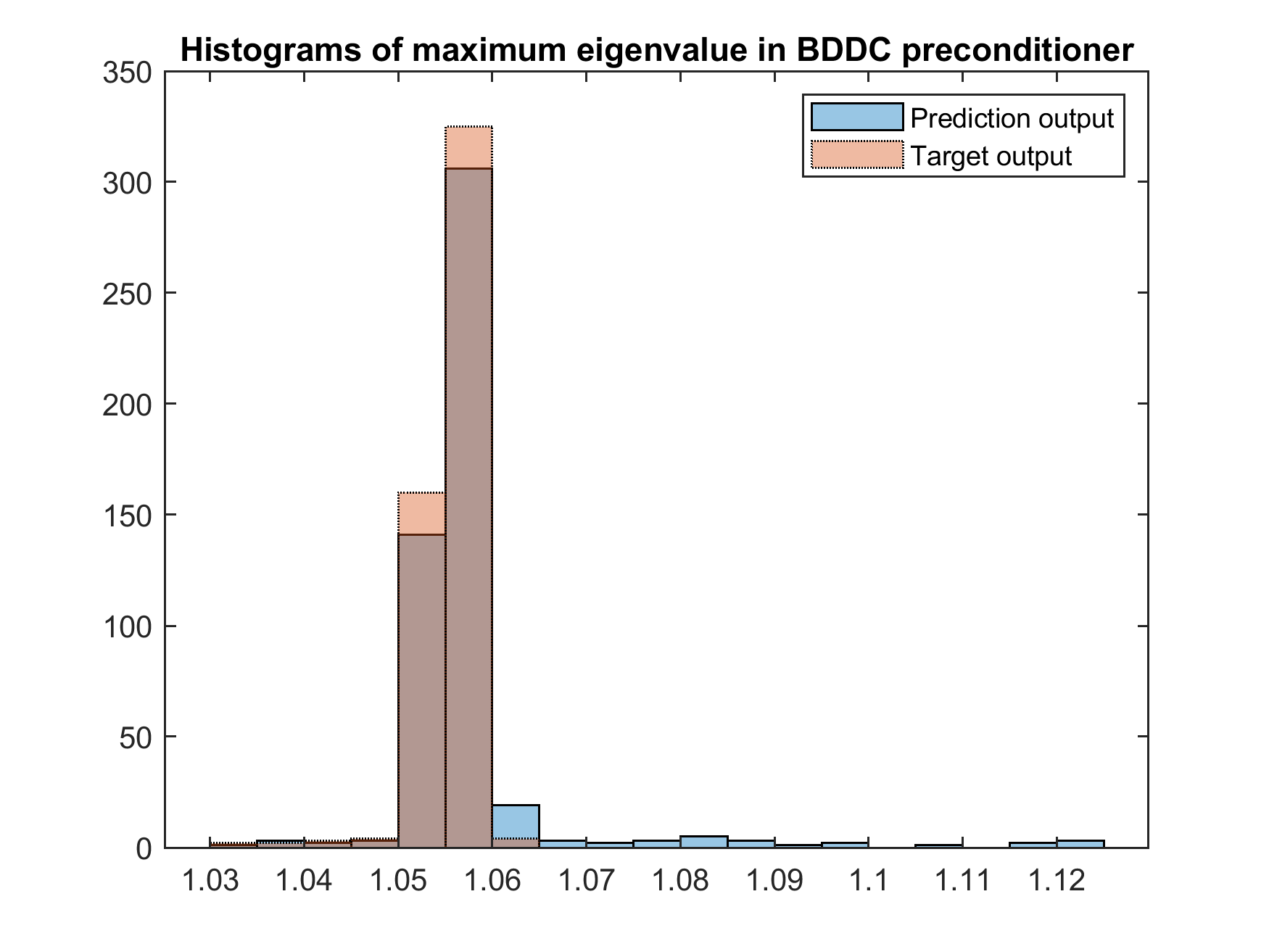}
\end{subfigure}
\caption{Comparisons in \rr{the performance of the} BDDC preconditioner  when expected functions of SPE10 Layer 35 (left column) and SPE10  Layer \rr{$35^*$} (right column) are used}
\label{fig:35compare}
\end{figure}

We first focus on the testing results of Layer 35. From the corresponding histograms in the left column of Figure \ref{fig:35compare}, besides the iteration number required for predicted eigenvectors are still usually larger by 1, we can observe a large area of overlapping, in particular, there is only about 10\% of data that are not overlapped in the histogram of maximum eigenvalues. Moreover, the small NRMSE and sMAPE again confirm with the graphical results. For generalization tests, we then consider two other testing sets with different stochastic behaviour and expected function.

\begin{figure}[h!]
\vspace{1cm}
\begin{subfigure}{.5\textwidth}
  \centering
  \includegraphics[width=1\linewidth]{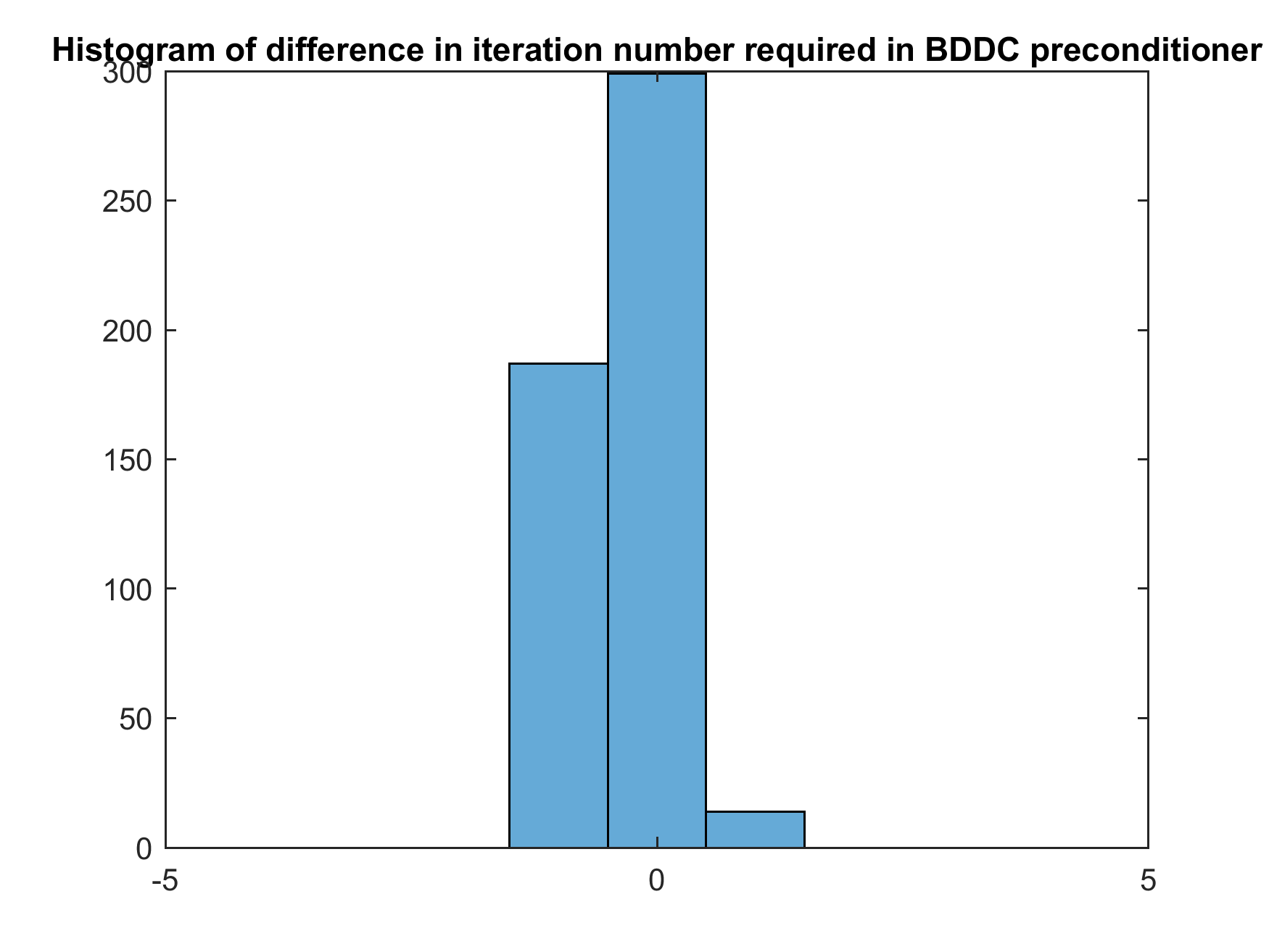}
\end{subfigure}%
\begin{subfigure}{.5\textwidth}
  \centering
  \includegraphics[width=1\linewidth]{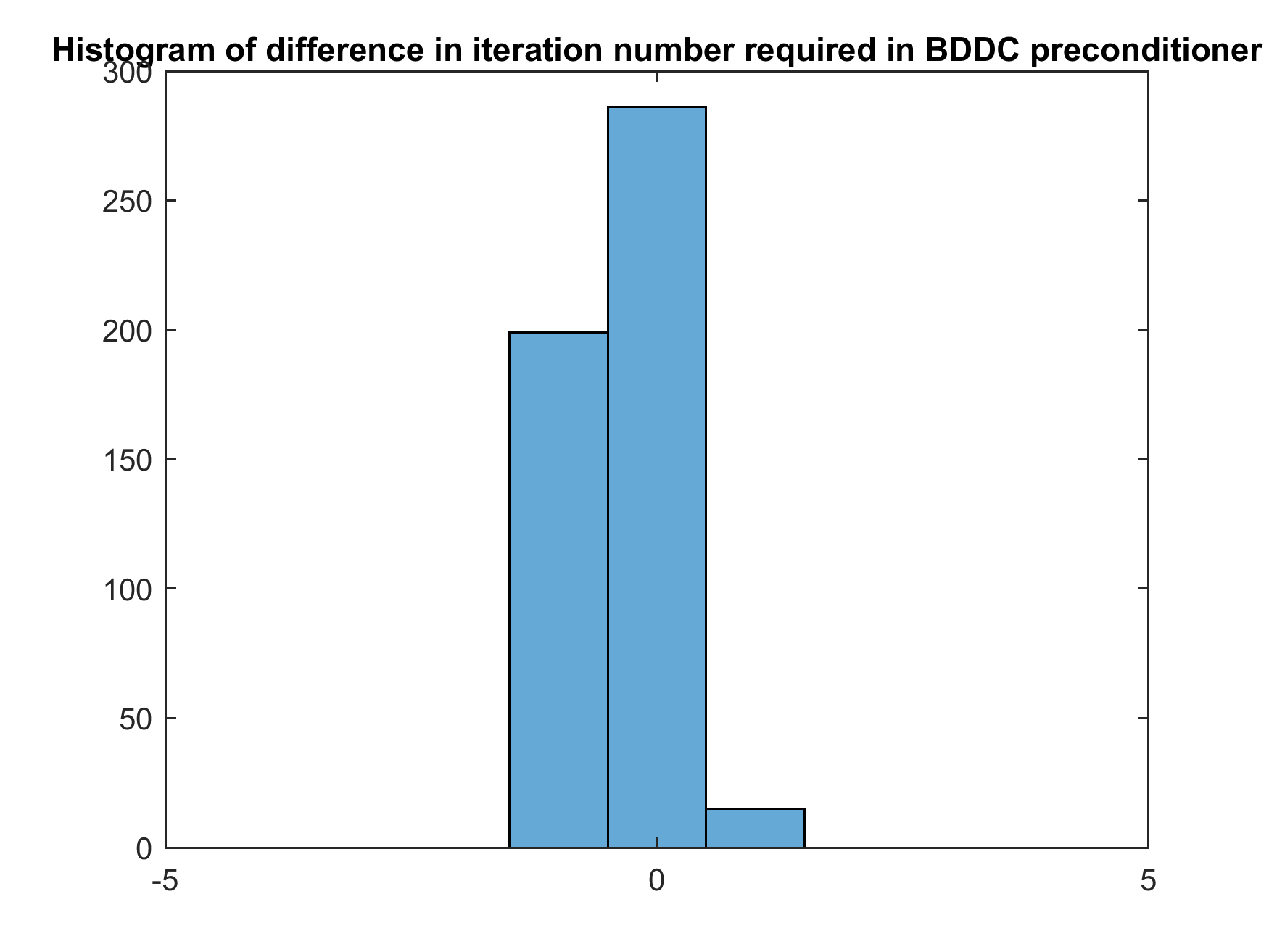}
\end{subfigure}

\begin{subfigure}{.5\textwidth}
  \centering
  \includegraphics[width=1\linewidth]{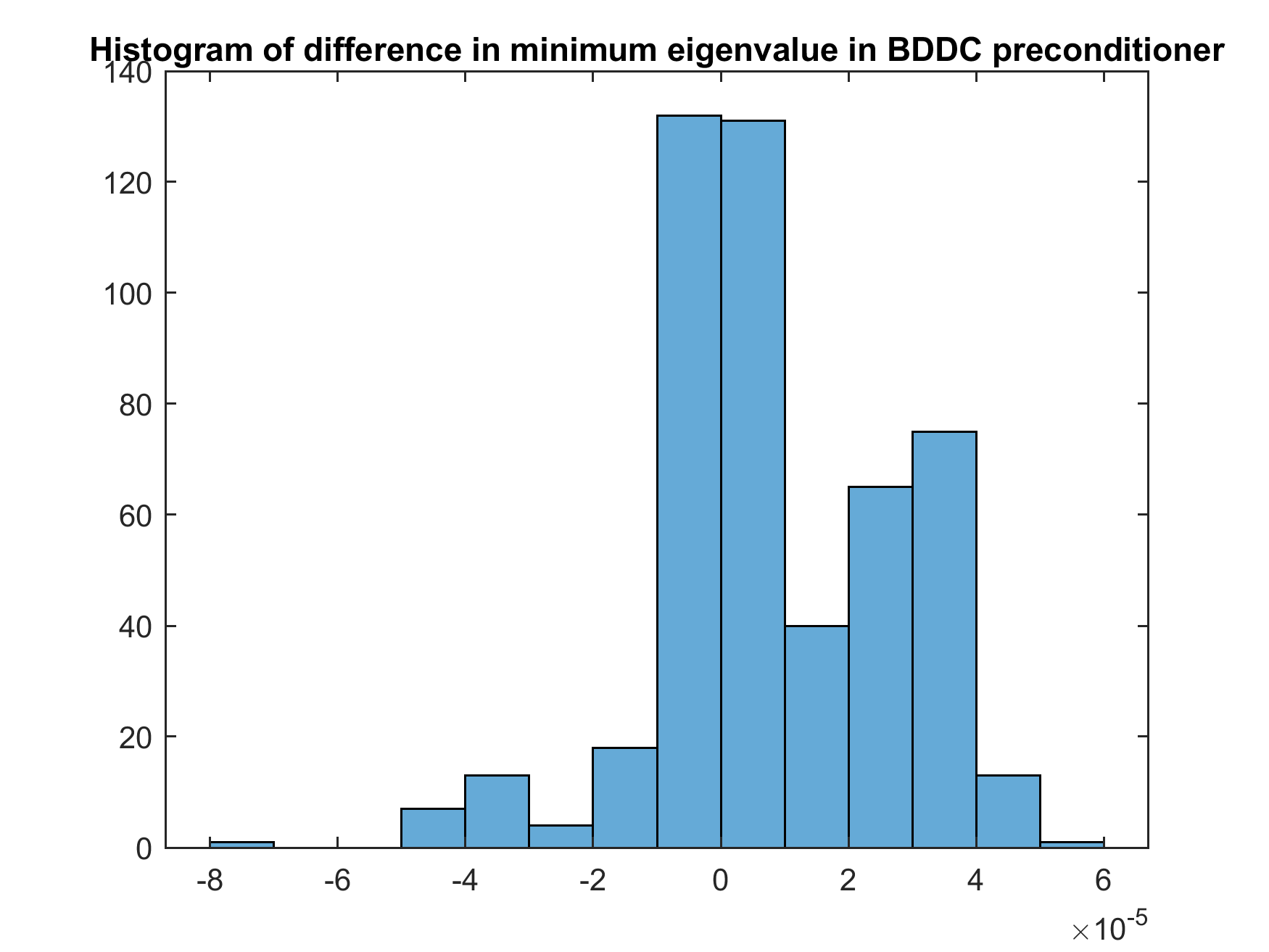}
\end{subfigure}%
\begin{subfigure}{.5\textwidth}
  \centering
  \includegraphics[width=1\linewidth]{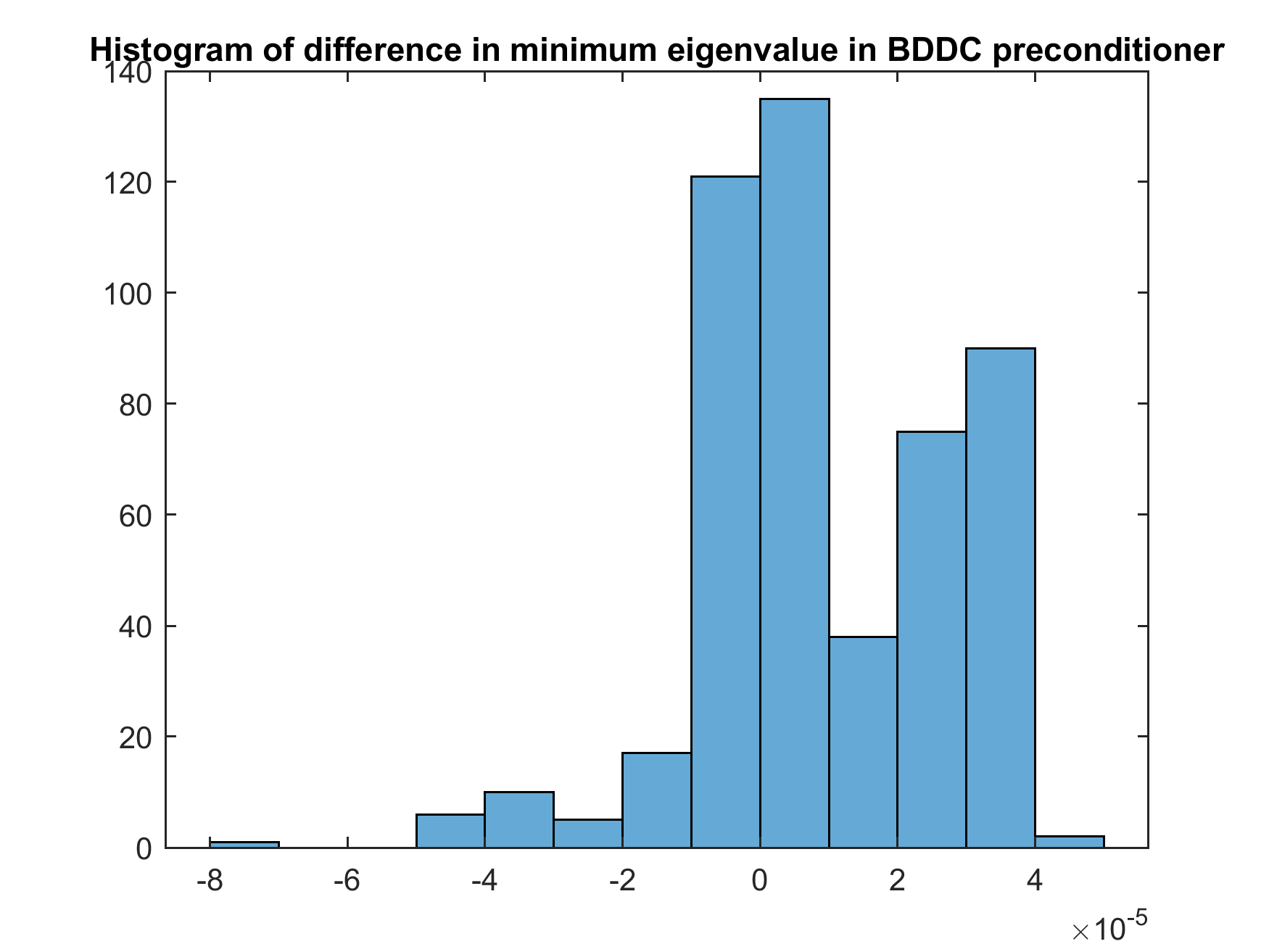}
\end{subfigure}

\begin{subfigure}{.5\textwidth}
  \centering
  \includegraphics[width=1\linewidth]{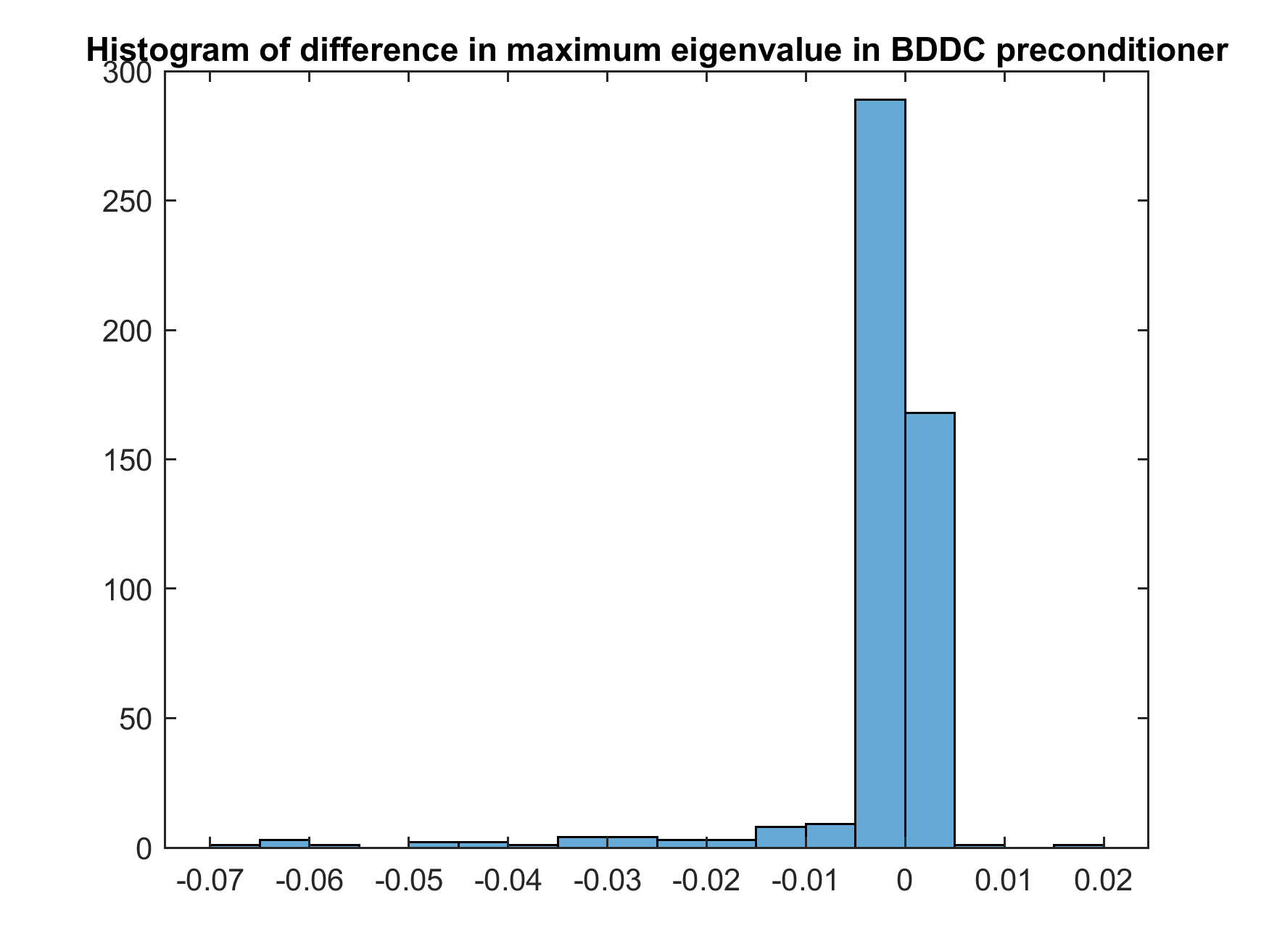}
\end{subfigure}%
\begin{subfigure}{.5\textwidth}
  \centering
  \includegraphics[width=1\linewidth]{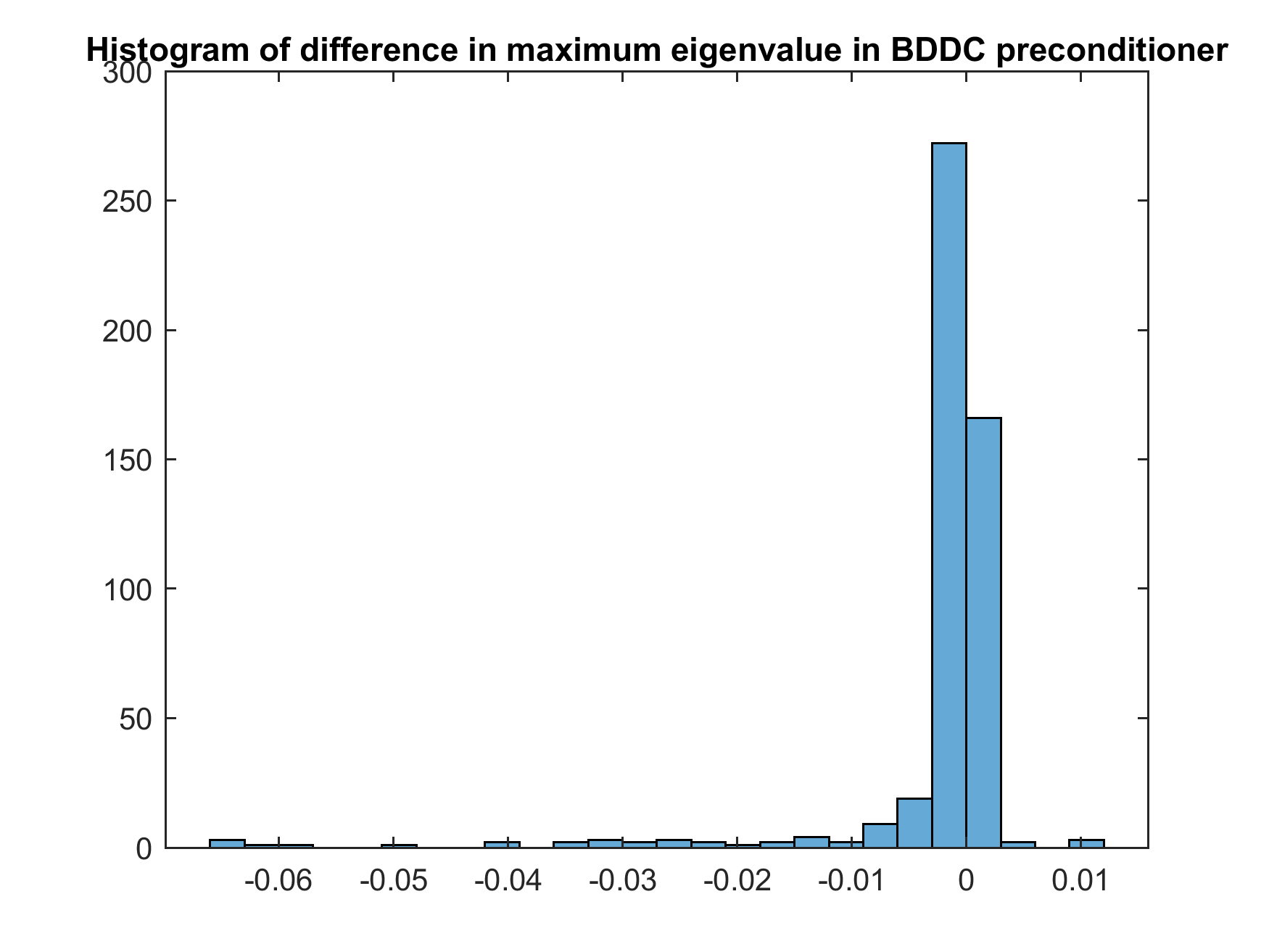}
\end{subfigure}
\caption{Differences in \rr{the performance of the} BDDC preconditioner when expected functions of SPE10 Layer 35 (left column) and SPE10  Layer \rr{$35^*$} (right column) are used}
\label{fig:35diff}
\end{figure}

When the stochastic behaviour of Layer 35 is changed to $(\sigma_K^2,\eta_1,\eta_2)=(1,0.2,0.125)$, we list the corresponding sMAPE and $l_\infty$ norms of differences in iteration numbers, maximum eigenvalues and minimum eigenvalues of preconditioners in the second row of Table \ref{test:Layer35}. We can see the values of errors are very similar to the first row and the NRMSE is even smaller. Therefore, we expect that the Layer 35 neural network can generalize the coefficient functions with similar stochastic properties, and give a good approximation on dominant \rr{eigenvectors for the BDDC preconditioner}. It is then verified by the histograms in the right column of Figure \ref{fig:35compare}. Moreover, a large population of samples are concentrated around 0 in the histograms of difference between \rr{the} estimated preconditioner and target preconditioner in Figure \ref{fig:35diff}. All these show our method performs well with stochastic oscillatory and high contrast coefficients, and also coefficient functions with similar stochastic properties.

\begin{figure}[h!]
\begin{subfigure}{.5\textwidth}
  \centering
  \includegraphics[width=1\linewidth]{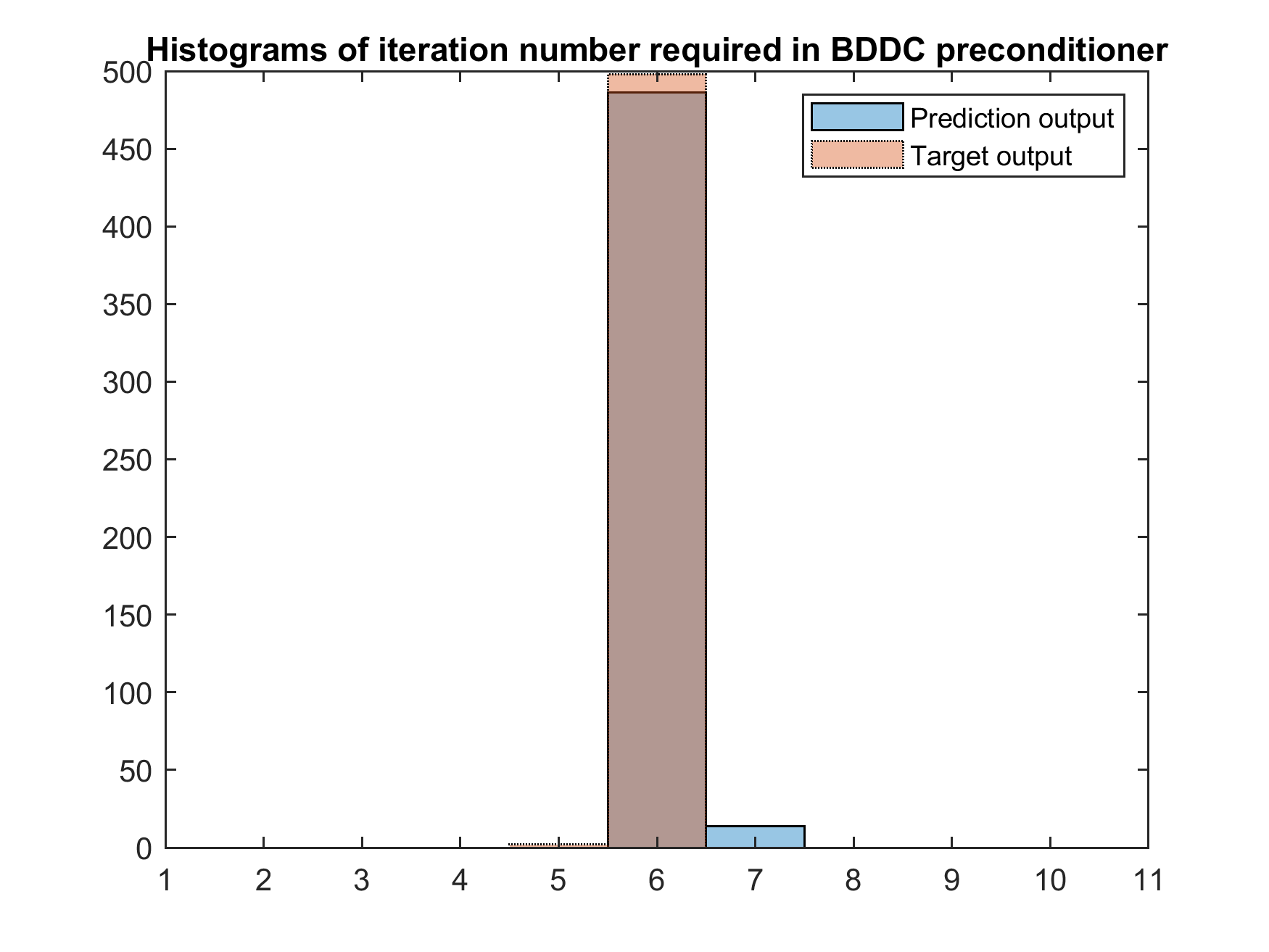}
\end{subfigure}%
\begin{subfigure}{.5\textwidth}
  \centering
  \includegraphics[width=1\linewidth]{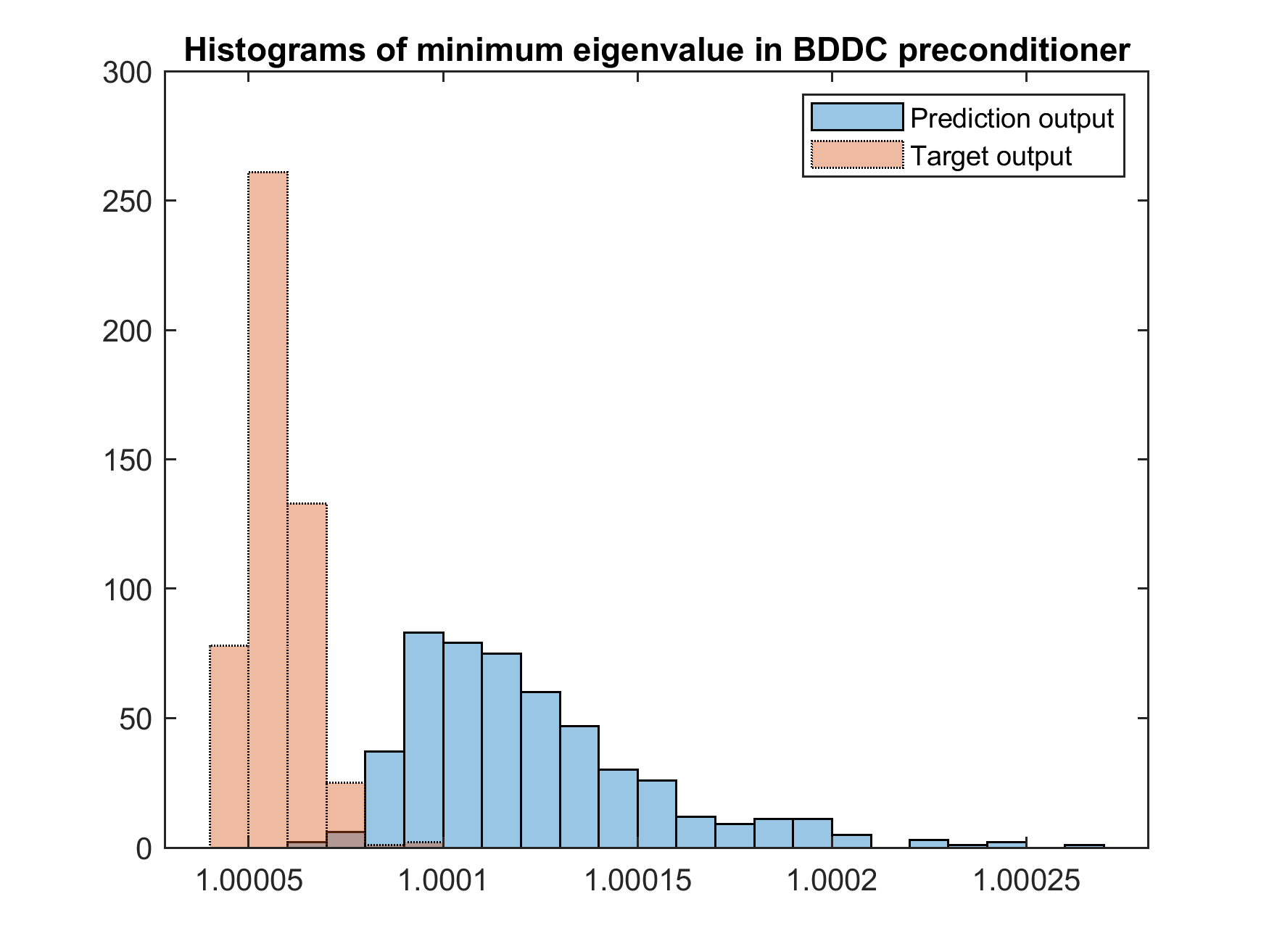}
\end{subfigure}

\hspace*{4cm}
\begin{subfigure}{.5\textwidth}
  \centering
  \includegraphics[width=1\linewidth]{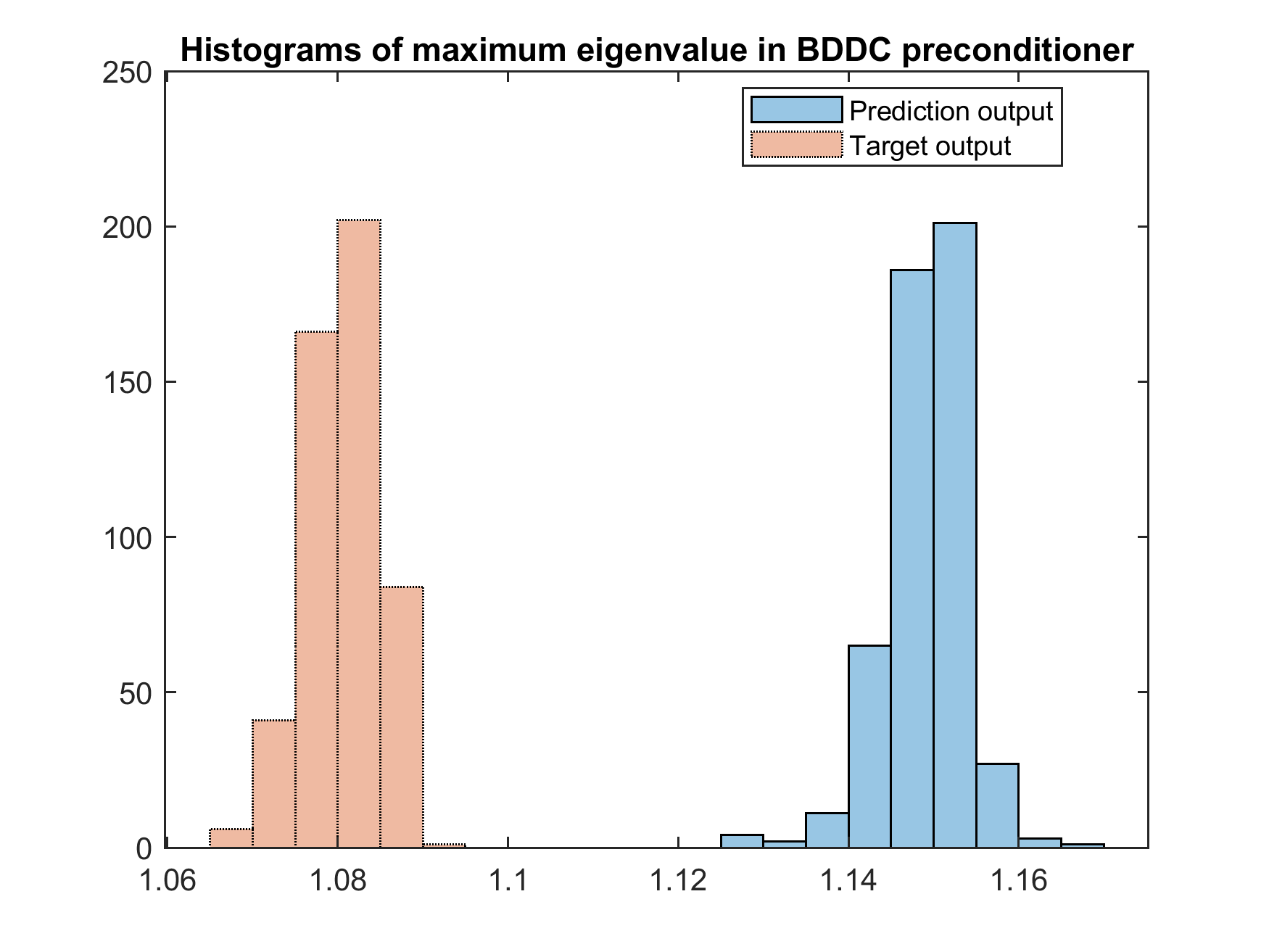}
\end{subfigure}%
\caption{Comparison in \rr{the performance of the} BDDC preconditioner  when expected function of SPE10 Layer 34 is used}
\label{fig:34compare}
\end{figure}

Finally, we compare the performances when the expected function is changed to Layer 34. Note, the NRMSE, sMAPE and $l_\infty$ errors of Layer 34 all are much higher than \rr{those} of Layer 35, which are crucial clues for \rr{a worse} performance. Although the permeability fields \rr{of} Layer 34 and Layer 35 share some common properties and have a near magnitude, unlike previous results for Layer 35 and Layer \rr{$35^*$}, we can see in Figure \ref{fig:34compare}, the prediction output and target output on maximum and minimum eigenvalues of preconditioners are even two populations with different centres and only \rr{a} little overlapping area. The main reason is a minor difference in $K(\boldsymbol{x},\omega)$ can already cause a huge change after taking exponential function. Therefore, when $E[K](\boldsymbol{x})$ is changed to another entirely different mean permeability with just some common properties, the feedforward neural network may not be able to capture these new characteristics, which will be part of our future research.

\clearpage

\section{\label{sec:conclusion}Conclusion}

A new learning adaptive BDDC algorithm is introduced and it consists of three main parts, which are the Gaussian random variables in Karhunen-Lo\`eve expansion, artificial neural network and the dominant eigenvectors obtained from the adaptive BDDC algorithm in the coarse spaces. The considered neural network acts as a fast computational tool, which turns the input Gaussian random variables into predicted dominant eigenvectors as output. \rr{In addition,} the neural network in the proposed algorithm is suitable for other permeability coefficients with similar stochastic properties for generalization purpose. Numerical results confirm the efficiency of the proposed algorithm and the generalization abilities. On the other hand, we currently just use the simplest feedforward neural network structure without \rr{retraining} when new characteristics are entered. In our future plan, we will improve the neural network and include deep learning technique for higher accuracy in prediction and preconditioning, and also better generalization ability on wider aspects of applications.


\bibliographystyle{plain}
\bibliography{adaptive-ref}

\begin{thebibliography}{10}

\bibitem{Parallel:Sum}
W.~N. Anderson, Jr. and R.~J. Duffin.
\newblock Series and parallel addition of matrices.
\newblock {\em J. Math. Anal. Appl.}, 26:576--594, 1969.

\bibitem{babuvska2007stochastic}
Ivo Babu{\v{s}}ka, Fabio Nobile, and Ra{\'u}l Tempone.
\newblock A stochastic collocation method for elliptic partial differential
  equations with random input data.
\newblock {\em SIAM Journal on Numerical Analysis}, 45(3):1005--1034, 2007.

\bibitem{babuska2004galerkin}
Ivo Babuska, Ra{\'u}l Tempone, and Georgios~E Zouraris.
\newblock Galerkin finite element approximations of stochastic elliptic partial
  differential equations.
\newblock {\em SIAM Journal on Numerical Analysis}, 42(2):800--825, 2004.

\bibitem{brunton2020machine}
Steven~L Brunton, Bernd~R Noack, and Petros Koumoutsakos.
\newblock Machine learning for fluid mechanics.
\newblock {\em Annual Review of Fluid Mechanics}, 52:477--508, 2020.

\bibitem{chung2021multi}
Eric Chung, Wing~Tat Leung, Sai-Mang Pun, and Zecheng Zhang.
\newblock A multi-stage deep learning based algorithm for multiscale model
  reduction.
\newblock {\em Journal of Computational and Applied Mathematics}, page 113506,
  2021.

\bibitem{chung2018adaptive}
Eric~T Chung, Yalchin Efendiev, and Wing~Tat Leung.
\newblock An adaptive generalized multiscale discontinuous galerkin method for
  high-contrast flow problems.
\newblock {\em Multiscale Modeling \& Simulation}, 16(3):1227--1257, 2018.

\bibitem{chung2014adaptive}
Eric~T Chung, Yalchin Efendiev, and Guanglian Li.
\newblock An adaptive gmsfem for high-contrast flow problems.
\newblock {\em Journal of Computational Physics}, 273:54--76, 2014.

\bibitem{BDDC-D-03}
Clark~R. Dohrmann.
\newblock A preconditioner for substructuring based on constrained energy
  minimization.
\newblock {\em SIAM J. Sci. Comput.}, 25(1):246--258, 2003.

\bibitem{Pechstein:slide:2013}
Clark~R. Dohrmann and Clemens Pechstein.
\newblock Modern domain decomposition solvers: {BDDC}, deluxe scaling, and an
  algebraic approach(2013),
  http://people.ricam.oeaw.ac.at/c.pechstein/pechstein-bddc2013.pdf.

\bibitem{dostert2006coarse}
Paul Dostert, Yalchin Efendiev, Thomas~Y Hou, and Wuan Luo.
\newblock Coarse-gradient langevin algorithms for dynamic data integration and
  uncertainty quantification.
\newblock {\em Journal of computational physics}, 217(1):123--142, 2006.

\bibitem{ghanem2003stochastic}
Roger~G Ghanem and Pol~D Spanos.
\newblock {\em Stochastic finite elements: a spectral approach}.
\newblock Courier Corporation, 2003.

\bibitem{heinlein2021combining}
Alexander Heinlein, Axel Klawonn, Martin Lanser, and Janine Weber.
\newblock Combining machine learning and domain decomposition methods for the
  solution of partial differential equations - a review.
\newblock {\em GAMM-Mitteilungen}, 44(1):e202100001, 2021.

\bibitem{kim2017bddc}
Hyea~Hyun Kim, Eric Chung, and Junxian Wang.
\newblock {BDDC} and {FETI-DP} preconditioners with adaptive coarse spaces for
  three-dimensional elliptic problems with oscillatory and high contrast
  coefficients.
\newblock {\em Journal of Computational Physics}, 349:191--214, 2017.

\bibitem{kim2015bddc}
Hyea~Hyun Kim and Eric~T Chung.
\newblock A {BDDC} algorithm with enriched coarse spaces for two-dimensional
  elliptic problems with oscillatory and high contrast coefficients.
\newblock {\em Multiscale Modeling \& Simulation}, 13(2):571--593, 2015.

\bibitem{Klawonn:PAMM:2014}
Axel Klawonn, Patrick Radtke, and Oliver Rheinbach.
\newblock {FETI}-{DP} with different scalings for adaptive coarse spaces.
\newblock {\em Proceedings in Applied Mathematics and Mechanics}, 2014.

\bibitem{KRR-2015}
Axel Klawonn, Patrick Radtke, and Oliver Rheinbach.
\newblock A comparison of adaptive coarse spaces for iterative substructuring
  in two dimensions.
\newblock {\em Electron. Trans. Numer. Anal.}, 45:75--106, 2016.

\bibitem{kutz2017deep}
J~Nathan Kutz.
\newblock Deep learning in fluid dynamics.
\newblock {\em Journal of Fluid Mechanics}, 814:1--4, 2017.

\bibitem{LW-FETIDP-BDDC}
Jing Li and Olof~B. Widlund.
\newblock F{ETI}-{DP}, {BDDC}, and block {C}holesky methods.
\newblock {\em Internat. J. Numer. Methods Engrg.}, 66(2):250--271, 2006.

\bibitem{lu2021deepxde}
Lu~Lu, Xuhui Meng, Zhiping Mao, and George~Em Karniadakis.
\newblock Deepxde: A deep learning library for solving differential equations.
\newblock {\em SIAM Review}, 63(1):208--228, 2021.

\bibitem{makridakis1993accuracy}
Spyros Makridakis.
\newblock Accuracy measures: theoretical and practical concerns.
\newblock {\em International journal of forecasting}, 9(4):527--529, 1993.

\bibitem{BDDC-Mandel-Dohrmann-Tezaur}
Jan Mandel, Clark~R. Dohrmann, and Radek Tezaur.
\newblock An algebraic theory for primal and dual substructuring methods by
  constraints.
\newblock {\em Appl. Numer. Math.}, 54(2):167--193, 2005.

\bibitem{moller1993scaled}
Martin~Fodslette M{\o}ller.
\newblock A scaled conjugate gradient algorithm for fast supervised learning.
\newblock {\em Neural Networks}, 6(4):525--533, 1993.

\bibitem{schwab2006karhunen}
Christoph Schwab and Radu~Alexandru Todor.
\newblock Karhunen-{L}o{\`e}ve approximation of random fields by generalized
  fast multipole methods.
\newblock {\em Journal of Computational Physics}, 217(1):100--122, 2006.

\bibitem{TW-Book}
Andrea Toselli and Olof Widlund.
\newblock {\em Domain decomposition methods---algorithms and theory}, volume~34
  of {\em Springer Series in Computational Mathematics}.
\newblock Springer-Verlag, Berlin, 2005.

\bibitem{vasilyeva2020learning}
Maria Vasilyeva, Wing~T Leung, Eric~T Chung, Yalchin Efendiev, and Mary
  Wheeler.
\newblock Learning macroscopic parameters in nonlinear multiscale simulations
  using nonlocal multicontinua upscaling techniques.
\newblock {\em Journal of Computational Physics}, 412:109323, 2020.

\bibitem{wang2008karhunen}
Limin Wang.
\newblock {\em Karhunen-Lo{\`e}ve expansions and their applications.}
\newblock PhD thesis, London School of Economics and Political Science (United
  Kingdom), 2008.

\bibitem{wang2020deep}
Yating Wang, Siu~Wun Cheung, Eric~T Chung, Yalchin Efendiev, and Min Wang.
\newblock Deep multiscale model learning.
\newblock {\em Journal of Computational Physics}, 406:109071, 2020.

\bibitem{wheeler2011multiscale}
Mary~F Wheeler, Tim Wildey, and Ivan Yotov.
\newblock A multiscale preconditioner for stochastic mortar mixed finite
  elements.
\newblock {\em Computer methods in applied mechanics and engineering},
  200(9-12):1251--1262, 2011.

\bibitem{yeung2020deep}
Tak Shing~Au Yeung, Eric~T Chung, and Simon See.
\newblock A deep learning based nonlinear upscaling method for transport
  equations.
\newblock {\em arXiv preprint arXiv:2007.03432}, 2020.

\bibitem{zhang2004efficient}
Dongxiao Zhang and Zhiming Lu.
\newblock An efficient, high-order perturbation approach for flow in random
  porous media via {K}arhunen-{L}o{\`e}ve and polynomial expansions.
\newblock {\em Journal of Computational Physics}, 194(2):773--794, 2004.

\bibitem{zhao2016spectral}
Wenzhi Zhao and Shihong Du.
\newblock Spectral--spatial feature extraction for hyperspectral image
  classification: A dimension reduction and deep learning approach.
\newblock {\em IEEE Transactions on Geoscience and Remote Sensing},
  54(8):4544--4554, 2016.

\end{thebibliography}

\end{document}